\documentclass[11pt, oneside]{article}
\usepackage{amsfonts}
\usepackage{mathrsfs}
\usepackage{float}
\usepackage{url}
\usepackage{graphicx}
\usepackage{epstopdf}
\usepackage{caption}
\usepackage{array}
\usepackage{amsmath}
\usepackage{bm}
\usepackage{booktabs}
\usepackage{multirow}
\usepackage{comment}
\usepackage[table]{xcolor}
\definecolor{wg}{HTML}{1AAD19}

\usepackage{latexsym,amssymb,amsmath,amstext}
\usepackage{mathtools}
\usepackage{amsthm}
\usepackage{indentfirst}
\usepackage{color}
\usepackage{tabularx}
\usepackage{subfigure}
\usepackage[ruled,vlined,linesnumbered]{algorithm2e}
\usepackage{enumitem}
\usepackage[colorlinks,linkcolor=blue,anchorcolor=blue,citecolor=blue]{hyperref}
\usepackage[authoryear,round,comma,sort&compress]{natbib}

\numberwithin{equation}{section}
\allowdisplaybreaks

\newtheorem{theorem}{\color{black}\indent Theorem}[section]

\newtheorem{definition}{\color{black}\indent Definition}[section]
\newtheorem{proposition}{\color{black}\indent Proposition}[section]

\makeatletter
\def\th@remark{%
  \normalfont
  \let\@begintheorem\@afterindenttrue\@afterheading\flushleft
  \let\@endtheorem\endflushleft
}
\makeatother

\theoremstyle{plain}

\newtheorem{corollary}{\color{black}\indent Corollary}[section]
\newtheorem{example}{\color{black}\indent Example}[section]

\begin{document}

\title{Additive Matrix Integer-Valued Autoregressive Model}
\author{\\Kaiyan Cui$^{1}$, Yikai Hu$^{1}$, Tianyun Guo$^{1}$\\
{$^{1}$\small School of Mathematics and Statistics, Shanxi University, Taiyuan 030006, China.}\\
{\small (E-mail: cuikaiyan@sxu.edu.cn )}
}

\date{}
\maketitle
\noindent {\bf Abstract}:
Contemporary data-driven and technology-integrated era, various matrix-valued integer-valued time series, such as cross-regional crime statistics, multi-category sales records, and network traffic matrices, exhibit high dimensionality, complex structures, and strong row-column intertwined dependencies. Although the existing matrix integer-valued autoregressive (MINAR) model provides a framework that directly handles matrix data and captures bidirectional row-column dependencies, it suffers from limited interpretability and inflexible structural representation, as its parameters often lack clear empirical meaning and the model cannot separately distinguish the effects arising from rows, columns, and lagged dynamics. To overcome these drawbacks, this paper proposes the additive matrix integer-valued autoregressive (Add-MINAR) model. By introducing an additive structure that explicitly decomposes the matrix response into row effects, column effects, and lagged effects, the proposed model not only preserves the matrix-valued nature but also significantly enhances parameter interpretability and structural flexibility. Two estimation methods, namely projection estimation and iterative conditional least squares estimation, are developed for parameter identification and inference, and their asymptotic properties, including consistency and asymptotic normality, are rigorously established. Simulation results show that the iterative conditional least squares estimator generally outperforms the projection estimator in most scenarios. Empirical analysis of Chicago crime data further demonstrates that the Add-MINAR model achieves superior in-sample fitting and out-of-sample forecasting performance compared to benchmark models such as MINAR, making it particularly suitable for practical applications with explicit row-column interaction features.

\vspace{0.3cm}

\thispagestyle{empty}
\section{Introduction}

Integer-valued time series are ubiquitous across various scientific and socio-economic disciplines, including epidemiology, economics, engineering, and public safety. Empirical examples include the spatial-temporal distribution of traffic accidents, daily retail sales volumes across different product categories, and hospital admission counts. Standard continuous time series models, such as the autoregressive moving average (ARMA) framework, are theoretically inadequate for count data. Applying continuous models to integer variables often yields nonsensical fractional forecasts and fails to account for the intrinsic distributional properties of count data, such as overdispersion, volatility clustering, or an excess of zeros. To address these fundamental issues, \citet{Steutel1979} introduced the binomial thinning operator, providing a probabilistic analogue to scalar multiplication for integer variables. Building upon this foundation, \citet{Al-Osh1987} proposed the first-order integer-valued autoregressive (INAR(1)) process. Over the subsequent decades, the INAR framework has been extensively generalized. Notable advancements include higher-order temporal dependence models \citep{Alzaid1990, Li1991}, alternative marginal specifications such as negative binomial or geometric distributions to accommodate overdispersion \citep{Al-Osh1992, Ristic2009}, frameworks handling zero-inflated phenomena \citep{Jazi2012}, and integer-valued generalized autoregressive conditional heteroskedasticity processes designed for volatility clustering \citep{Ferland2006, Zhu2011}. 

Meanwhile, the rise of complex and structured data, such as multi-regional crime statistics or cross-category retail sales matrices, has exposed critical limitations of univariate count models, which fail to account for inherent cross-sectional correlations present in the data. To address this gap, \citet{Franke1993} initiated the multivariate extension, which was subsequently formalized into the multivariate generalized integer-valued autoregressive (MGINAR) model by \citet{Latour1997}. Later research thoroughly investigated the dependency structures and practical applications of these multivariate systems \citep{Pedeli2011, Pedeli2013a}, while advancing estimation techniques such as composite likelihood methods to mitigate associated computational burdens \citep{Pedeli2013b, Zhang2017}. However, when applied to matrix integer-valued time series, these traditional multivariate approaches rely on vectorizing the observation matrix into a single long column vector. This vectorization procedure inherently destroys the underlying spatial and cross-sectional topologies of the data. More critically, it induces a severe parameter proliferation that renders estimation computationally intractable, significantly increases the risk of overfitting, and proves statistically inefficient in complex empirical applications.


The theoretical and computational limitations of vector-based models have driven researchers to shift focus toward methodologies that preserve the intrinsic matrix structure. In the continuous domain, matrix-valued factor models \citep{Wang2019} and continuous matrix autoregressive (MAR) models \citep{Chen2021} have successfully characterized row and column dynamics parsimoniously, avoiding the parameter explosion associated with vector autoregressions. Inspired by this structural paradigm, the matrix integer-valued autoregressive (MINAR) model was recently introduced for count data by \citet{Xu2025,Cui2026}. By employing a novel matrix thinning operator, the MINAR framework directly processes matrix observations, thereby significantly reducing the parameter space and retaining structural integrity. Nevertheless, existing MINAR models depend strictly on a coupled bilinear product specification. This structural formulation mathematically conflates the row and column effects into a single interactive process, thereby severely restricting the empirical interpretability of the parameters. For instance, in regional crime analysis, a bilinear model cannot explicitly isolate whether a fluctuation in a specific spatial-temporal entry is driven by cross-regional spillover (a column effect), cross-crime contagion (a row effect), or the intrinsic historical inertia of that specific variable (a lagged effect).

To address these critical structural limitations of the bilinear MINAR framework and advance the modeling of matrix-valued count data, this paper proposes the additive matrix integer-valued autoregressive (Add-MINAR) model. By incorporating an innovative additive specification, the proposed model explicitly decomposes the matrix response into three independent dynamic components: row interactive effects (matrix $\mathbf{A}$), column interactive effects (matrix $\mathbf{B}$), and idiosyncratic autoregressive lagged effects (matrix $\mathbf{C}$). To theoretically ensure parameter identifiability and prevent estimation ambiguity within this additive architecture, we impose strict zero-diagonal constraints on the interactive coefficient matrices $\mathbf{A}$ and $\mathbf{B}$. This structural decomposition not only inherits the parameter reduction advantages of matrix modeling but also fundamentally enhances the interpretability of the parameters. By separating these effects, the model allows researchers to test specific hypotheses about spatial and cross-sectional dynamics independently, allowing practitioners to transparently trace the distinct dynamic drivers of the modeled system.

The primary contributions of this paper are threefold. First, we establish the probabilistic foundations of the Add-MINAR model utilizing element-wise sparsification operators. We rigorously derive its stationarity and causality conditions, ensuring that the additive structure maintains the theoretical soundness of traditional integer-valued processes. Second, we develop two distinct parameter estimation methodologies: projection estimation (PROJ) and iterative conditional least squares estimation (ICLSE). By formulating the zero-diagonal constraints as linear equality restrictions, the ICLSE is systematically solved via the Lagrangian multiplier method, iteratively optimizing the objective function. We systematically establish the consistency and asymptotic normality of both estimators, providing a solid theoretical guarantee for statistical inference. Third, comprehensive Monte Carlo simulations and an empirical application to Chicago crime statistics demonstrate the practical superiority of the proposed framework. The results indicate that the Add-MINAR model, particularly under the ICLSE procedure, exhibits robust parameter estimation, effectively prevents overfitting, and consistently outperforms traditional MGINAR and bilinear MINAR benchmarks in both in-sample fitting accuracy and out-of-sample forecasting stability.

The remainder of this paper is organized as follows. Section 2 formally defines the Add-MINAR model and discusses its probabilistic properties. Section 3 delineates the PROJ and ICLSE methodologies in detail. Section 4 establishes the asymptotic theory for the proposed estimators. Section 5 presents finite-sample simulation results and the empirical crime data analysis. Finally, Section 6 concludes the paper. All technical proofs are relegated to the Appendix.

\vspace{0.3cm}

\noindent{\it Notation.} 
We use $\otimes$ to denote the Kronecker product, $\|\cdot\|_F^{2}$ the Frobenius norm of a matrix, and $\|\cdot\|_2$ the spectral norm of a matrix. The bold capital letters, e.g., $\mathbf{A}$, denote matrices, and the bold lowercase letters, e.g., $\mathbf{v}$, denote vectors. $a_{i,j}$ denotes the element in the \(i\)-th row and \(j\)-th column of matrix $\mathbf{A}$. Let $\operatorname{vec}(\cdot)$ be the vectorization of a matrix by stacking its columns. Define $\rho(\cdot)$ as the spectral radius of a matrix, i.e., the maximum modulus of its eigenvalues. For any vector $\mathbf{v} \in \mathbb{R}^{m}$, $\operatorname{diag}(\mathbf{v}) \in \mathbb{R}^{m \times m}$ denotes the diagonal matrix formed by the elements of $\mathbf{v}$. For any matrix $\mathbf{A} \in \mathbb{R}^{m \times m}$, $\operatorname{diag}(\mathbf{A}) \in \mathbb{R}^{m}$ denotes the vector formed by the diagonal elements of $\mathbf{A}$; while $\operatorname{Diag}(\mathbf{A}) \in \mathbb{R}^{m \times m}$ denotes the diagonal matrix formed by the diagonal elements of $\mathbf{A}$. The trace and rank of the matrix $\mathbf{A}$ are denoted by $\operatorname{tr}(\mathbf{A})$ and $\operatorname{rank}(\mathbf{A})$, respectively. We use $\odot$ and $\oslash$ to denote the element-wise multiplication and element-wise division of matrices, respectively. The $m$-dimensional identity matrix and the $(m, n)$-dimensional all-zero matrix are denoted by $\mathbf{I}_m$ and $\mathbf{0}_{m \times n}$, respectively. Let $\mathbf{e}_{d,k} \in \mathbb{R}^d$ denote the column vector with the $k$-th entry equal to $1$ and all other entries equal to $0$. Furthermore, $\mathcal{R}_{m,n}: \mathbb{R}^{mn} \to \mathbb{R}^{m \times n}$ denotes the reshaping operator that maps a vectorized input back into its corresponding $m \times n$ matrix format.

\section{Additive Matrix Integer-Valued Autoregressive Model}

In this section, we formally introduce the Add-MINAR framework. We first introduce the rigorous mathematical definitions of the matrix thinning operators utilized within the proposed model. Subsequently, we present the overarching model specifications and establish the fundamental matrix notation. Finally, we provide an in-depth discussion on three fundamental aspects of the framework: parameter identifiability, nested relationships among model variants, and the vectorized representation required for parameter estimation.

\subsection{Definitions of Thinning Operators}\label{subsec:operators}

To establish a rigorous foundation, we sequentially introduce these operators from the univariate case to the matrix dimension. 

\begin{definition}[Univariate General Thinning Operator]\label{def:binomial_thinning} 
Let $\{\xi_{k}^{(\alpha)}\}_{k \in \mathbb{N}}$ be a sequence of independent and identically distributed ({\rm i.i.d.}) non-negative integer-valued random variables with finite mean $\alpha \geq 0$ and finite variance, independent of $x$, an $\mathbb{N}_0$-valued random variable. The generalized Steutel and van Harn thinning operator, denoted by $\alpha\circ$, is defined as
\[
\alpha\circ x\coloneqq 
\begin{cases} 
\sum_{k=1}^{x}\xi_{k}^{(\alpha)}, & x>0, \\ 
0, & x=0. 
\end{cases}
\]
The sequence $\{\xi_{k}^{(\alpha)}\}_{k \in \mathbb{N}}$ is referred to as a counting sequence.
\end{definition}

In this paper, we specifically adopt the Poisson thinning operator, assuming $\xi_{k}^{(\alpha)} \sim \mathrm{Poisson}(\alpha)$ for the counting sequence in Definition \ref{def:binomial_thinning}.

\begin{definition}[Multivariate Vector Thinning] \label{def:vector_thinning}
For a $m\times m$ matrix $\mathbf{A}=(\alpha_{i,j})_{m\times m}\in \mathbb{R}_{\geq0}^{m\times m}$ and an $\mathbb{N}_0^{m}$-valued random vector $\mathbf{v}=(v_1,\ldots,v_m)^\top$, define the vector (multivariate) thinning operator $\ast$ by
\[\mathbf{A}\ast\mathbf{v}\coloneqq\left(\sum_{k=1}^{m}\alpha_{i,k}\circ v_{k}\right)_{1\leq i\leq m},\]
where the general thinning operator $\alpha_{i,k}\circ\cdot$ independently over $1\leq i,k\leq m$.
\end{definition}

\begin{definition}[Single Matrix Thinning] \label{def:matrix_thinning_single}
For a $m\times m$ matrix $\mathbf{A}=(\alpha_{i,j})_{m\times m}\in \mathbb{R}_{\geq0}^{m\times m}$, an $n\times n$ matrix $\mathbf{B}=(\beta_{i,j})_{n\times n}\in \mathbb{R}_{\geq0}^{n\times n}$ and an $\mathbb{N}_0^{m\times n}$-valued random matrix $\mathbf{Y}=(y_{i,j})_{m\times n}$, define the single matrix thinning operator $\circledast$ by
\[\mathbf{A}\circledast\mathbf{Y}\coloneqq \left(\sum_{k=1}^{m}\alpha_{i,k}\circ y_{k,j}\right)_{1\leq i\leq m,1\leq j\leq n},\]
\[\mathbf{Y}\circledast \mathbf{B}\coloneqq \left(\sum_{l=1}^{n}\beta_{l,j}\circ y_{i,l}\right)_{1\leq i\leq m,1\leq j\leq n}.\]
\end{definition}


\begin{definition}[Element-wise Sparsification]\label{def:elementwise_thinning}
For a $m\times n$ matrix $\mathbf{C}=(\alpha_{i,j})_{m\times n}\in\mathbb{R}_{\geq0}^{m\times n}$, an $\mathbb{N}_0^{m\times n}$-valued random matrix $\mathbf{Y}=(y_{i,j})_{m\times n}$, define the elementwise sparsification operator $\circledast~\cdot$ by
\[
(\mathbf{C}\circledast\cdot \mathbf{Y})_{i,j} \coloneqq \alpha_{i,j}\circ y_{i,j},
\]
where for $1\leq i\leq m$, $1\leq j\leq n$, the sparsification operators $\alpha_{i,j}\circ\cdot$ are independent of each other.
\end{definition}

\subsection{Model Specification}\label{sec:Model Specification}

The core motivation of the Add-MINAR framework is to explicitly decompose the dynamics of matrix-valued count data into row interactive effects, column interactive effects, and idiosyncratic lagged effects. To formalize this, let the parameter matrices be defined as:
\[
\begin{aligned}
\mathbf{A} &= (a_{i,j})_{1 \leq i, j \leq m} \in \mathbb{R}_{\geq 0}^{m \times m}, \quad \mathbf{B} = (b_{i,j})_{1 \leq i, j \leq n} \in \mathbb{R}_{\geq 0}^{n \times n}, \\
\mathbf{C} &= (c_{i,j})_{1 \leq i \leq m, 1 \leq j \leq n} \in \mathbb{R}_{\geq 0}^{m \times n}, \quad \mathbf{D} = (d_{i,j})_{1 \leq i \leq m, 1 \leq j \leq n} \in \mathbb{R}_{\geq 0}^{m \times n}.
\end{aligned}
\]
Let $\{\mathbf{X}_t\}_{t \in \mathbb{Z}}$ be the matrix integer-valued time series with elements $\mathbf{X}_t = (x_{t,i,j})_{m \times n}$, and $\{\mathbf{E}_t\}_{t \in \mathbb{Z}}$ be a sequence of independent and identically distributed non-negative integer-valued random matrices with elements $\mathbf{E}_t = (e_{t,i,j})_{m \times n}$, where $e_{t,i,j} \sim \mathrm{Poisson}(d_{i,j})$ for $t \in \mathbb{Z}$. Assume that for any $t_1 > t_2$, $\mathbf{E}_{t_1}$ is independent of $\mathbf{X}_{t_2}$. The three variants of the \textnormal{Add-MINAR(1)} model are defined as:
\begin{enumerate}
\item[\rm{(I)}] \(\displaystyle \mathbf{X}_t = \mathbf{A} \circledast \mathbf{X}_{t-1} + \mathbf{X}_{t-1} \circledast \mathbf{B}^\top + \mathbf{E}_t, \quad\) 
\textnormal{subject to } \(\operatorname{diag}(\mathbf{A}) = \mathbf{0}_{m \times 1};\)

\item[\rm{(II)}] \(\displaystyle \mathbf{X}_t = \mathbf{A} \circledast \mathbf{X}_{t-1} + \mathbf{X}_{t-1} \circledast \mathbf{B}^\top + \mathbf{E}_t, \quad\) 
\textnormal{subject to } \(\operatorname{diag}(\mathbf{B}) = \mathbf{0}_{n \times 1};\)

\item[\rm{(III)}] \(\displaystyle \mathbf{X}_t = \mathbf{A} \circledast \mathbf{X}_{t-1} + \mathbf{X}_{t-1} \circledast \mathbf{B}^\top + \mathbf{C} \circledast\cdot \mathbf{X}_{t-1} + \mathbf{E}_t, \) \\
\textnormal{subject to } \(\operatorname{diag}(\mathbf{A}) = \mathbf{0}_{m \times 1} \textnormal{ and } \operatorname{diag}(\mathbf{B}) = \mathbf{0}_{n \times 1}.\) 
\end{enumerate}

Extended from the Add-MAR model\citep{Zhang2024}, the proposed Add-MINAR framework, particularly Model (III), explicitly decouples the dynamics of the matrix-valued time series into three distinct and interpretable components. The coefficient matrix $\mathbf{A}$ serves as the row parameter matrix, capturing the cross-sectional inter-effects and spillovers across the rows of $\{\mathbf{X}_t\}_{t \in \mathbb{Z}}$. The matrix $\mathbf{B}$ acts as the column parameter matrix, reflecting the interactive effects across the columns. Meanwhile, the parameter matrix $\mathbf{C}$ accounts for the specific, independent lagged momentum (or historical inertia) of each individual cell.

The fundamental motivation for this additive decomposition becomes profoundly evident when contrasted with the standard MINAR model. Existing MINAR frameworks typically rely on a bilinear multiple matrix thinning structure, mathematically formulated as $\mathbf{A}\circledast\mathbf{X}_{t-1}\circledast \mathbf{B}^{\top}$. Under this bilinear specification, the expected historical dependency relies on a coupled double summation, where the dynamic of the $(i,j)$-th entry is driven by $\sum_{l=1}^{n}\sum_{k=1}^{m}(b_{j,l}a_{i,k})\circ x_{t-1,k,l}$. While mathematically compact, this formulation intrinsically intertwines the row effects ($a_{i,k}$) and column effects ($b_{j,l}$) into a single multiplicative interactive process. Consequently, in empirical applications, it becomes highly challenging to disentangle whether a specific variation is caused by row-wise contagion, column-wise spillover, or the purely autoregressive memory of the cell itself. By abandoning the bilinear product in favor of the additive specification ($\mathbf{A} \circledast \mathbf{X}_{t-1} + \mathbf{X}_{t-1} \circledast \mathbf{B}^\top + \mathbf{C} \circledast\cdot \mathbf{X}_{t-1}$), the Add-MINAR model successfully isolates these temporal and spatial dynamics, thereby fundamentally enhancing the structural flexibility and empirical interpretability of the parameters.

\subsection{Parameter Identifiability}\label{subsec:identifiability}

A critical feature of the proposed Add-MINAR framework is the imposition of zero-diagonal constraints: $\operatorname{diag}(\mathbf{A}) = \mathbf{0}_{m \times 1}$ and $\operatorname{diag}(\mathbf{B}) = \mathbf{0}_{n \times 1}$. These constraints are mathematically imperative to ensure the unique identifiability of the model parameters. 

To illustrate this necessity, consider the unconstrained element-wise expansion of the row and column interactive terms in Model (I). Based on the notation established in Model (I) specified in Section \ref{sec:Model Specification}, the expected value for the $(i, j)$-th entry of the generated matrix at time $t$ involves the following sum:
\[
x_{t,i,j} = \sum_{k=1}^m a_{i,k} \circ x_{t-1,k,j} + \sum_{s=1}^n b_{j,s} \circ x_{t-1,i,s} + e_{t,i,j}.
\]
When evaluating the specific case where $1 \leq k \leq m$ and $1 \leq s \leq n$ coincide at $k = i$ and $s = j$, the right-hand side yields a combined coefficient for the self-lagged term $x_{t-1, i, j}$:
\[
a_{i,i} \circ x_{t-1,i,j} + b_{j,j} \circ x_{t-1,i,j} \overset{\mathrm{d}}{=} (a_{i,i} + b_{j,j}) \circ x_{t-1,i,j}.
\]
Without structural restrictions, there exists an arbitrary scalar $\delta \neq 0$ such that:
\[
(a_{i,i} + \delta) + (b_{j,j} - \delta) = a_{i,i} + b_{j,j}.
\]
This arbitrary shifting of mass between $a_{i,i}$ and $b_{j,j}$ leads to an infinite number of parameter combinations yielding the exact same conditional expectation, thus rendering the estimation system non-identifiable. By forcing the diagonals of $\mathbf{A}$ and $\mathbf{B}$ to be zero ($\operatorname{diag}(\mathbf{A}) = \mathbf{0}_{m \times 1}$ and $\operatorname{diag}(\mathbf{B}) = \mathbf{0}_{n \times 1}$) in Models (I) and (II), respectively, and delegating the self-lagged dynamics entirely to the element-wise matrix $\mathbf{C}$ in Model (III), we completely eliminate this estimation ambiguity.

\subsection{Nested Relationships Among Model Variants}\label{subsec:relationships}

Analogous to the structural properties discussed in the continuous Add-MAR framework\citep{Zhang2024}, Models (I) and (II) mathematically serve as nested boundary cases of the generalized Model (III). Specifically, if the idiosyncratic lag matrix $\mathbf{C}$ exhibits a column-replicated structure (i.e., homogeneous across columns), its specific lag effects can be entirely absorbed by relaxing the zero-diagonal constraint on the row matrix $\mathbf{A}$. This absorption eliminates the independent matrix $\mathbf{C}$, thereby exactly degenerating Model (III) into Model (II). 

Symmetrically, if $\mathbf{C}$ assumes a row-replicated structure, its parameters can be proportionally absorbed into the unconstrained diagonal elements of the column matrix $\mathbf{B}$, reducing the framework to Model (I). Because Model (III) encompasses the most generalized and structurally complete parameterization, the theoretical derivations and empirical analyses in the remainder of this paper will focus exclusively on Model (III).

\subsection{Vectorized Representation}\label{subsec:vectorization}

To facilitate parameter estimation and establish asymptotic properties, we algebraically transform the matrix-variate Model (III) into a standard high-dimensional MGINAR(1) structure. Let $\operatorname{vec}(\cdot)$ denote the column-wise vectorization operator. Applying this operator to the observation system yields:
\begin{equation}\label{equ:Add-MINAR to MGINAR}
\operatorname{vec}(\mathbf{X}_t) = \mathbf{\Phi} \ast \operatorname{vec}(\mathbf{X}_{t-1}) +\operatorname{vec}(\mathbf{E}_t),
\end{equation}
where $\mathbf{\Phi} \in \mathbb{R}^{mn \times mn}$ functions as the global autoregressive coefficient matrix. 

Instead of explicitly partitioning $\mathbf{\Phi}$ into localized block matrices, we can compactly construct it by leveraging the algebraic properties of the matrix thinning operators. The overarching matrix $\mathbf{\Phi}$ is essentially the superposition of the row, column, and idiosyncratic effects, formulated elegantly via Kronecker products as:
\begin{equation}\label{equ:Phi_compact}
\mathbf{\Phi} = (\mathbf{I}_{n \times n} \otimes \mathbf{A}) + (\mathbf{B} \otimes \mathbf{I}_{m \times m}) + \operatorname{Diag}(\operatorname{vec}(\mathbf{C})),
\end{equation}
subject to the identification constraints $\operatorname{diag}(\mathbf{A}) = \mathbf{0}_{m \times 1}$ and $\operatorname{diag}(\mathbf{B}) = \mathbf{0}_{n \times 1}$. Here, $\operatorname{Diag}(\cdot)$ transforms a vector into a diagonal matrix. It is evident from \eqref{equ:Phi_compact} that the zero-diagonal restrictions on $\mathbf{A}$ and $\mathbf{B}$ ensure that the main diagonal of the global matrix $\mathbf{\Phi}$ is exclusively governed by $\mathbf{C}$.

To ensure the statistical stability of the Add-MINAR process, we impose the spectral radius condition $\rho(\mathbf{\Phi}) < 1$. This condition guarantees the existence of a unique, strictly stationary, and ergodic solution for $\{\mathbf{X}_t\}_{t \in \mathbb{Z}}$ \cite{Latour1997}, and further implies causality where the influence of historical innovations diminishes at a geometric rate.

\section{Parameter Estimation}

In this section, we proposes two estimation methodologies for the parameter matrices $\mathbf{A}$, $\mathbf{B}$, $\mathbf{C}$, and $\mathbf{D}$ within the Add-MINAR(1) model. The first is the projection estimation (PROJ) method, which directly estimates the coefficient matrix of the vectorized MGINAR(1) model and subsequently restores the structural matrix parameters through specialized projection operations. The second is the iterative conditional least squares estimation (ICLSE) method, which estimates the parameters through an alternating fixed-parameter approach and iterative optimization under equality constraints. To ensure clarity and prevent logical ambiguity, the estimates obtained via the PROJ method are denoted as \((\hat{\mathbf{A}}_{\mathrm{PROJ}}, \hat{\mathbf{B}}_{\mathrm{PROJ}}, \hat{\mathbf{C}}_{\mathrm{PROJ}}, \hat{\mathbf{D}}_{\mathrm{PROJ}})\), while those obtained via the ICLSE method are denoted as \((\hat{\mathbf{A}}, \hat{\mathbf{B}}, \hat{\mathbf{C}}, \hat{\mathbf{D}})\).

\subsection{Projection Estimation Method}\label{PROJ}

From \eqref{equ:Add-MINAR to MGINAR}, the Add-MINAR(1) model can be algebraically transformed into a standard MGINAR(1) model. The corresponding conditional expectation is expressed as:
\begin{align*}
\mathbb{E}(\operatorname{vec}(\mathbf{X}_t)\mid\operatorname{vec}(\mathbf{X}_{t-1}))
&=\mathbf{\Phi}\operatorname{vec}(\mathbf{X}_{t-1})+\operatorname{vec}(\mathbf{D})\\
&=(\mathbf{\Phi}\vdots\operatorname{vec}(\mathbf{D}))\binom{\operatorname{vec}(\mathbf{X}_{t-1})}{1}
=:\mathbf{\Psi}\operatorname{vec}(\mathbf{Z}_{t-1}),
\end{align*}
where
$$
\operatorname{vec}(\mathbf{Z}_{t-1})=\begin{pmatrix}\operatorname{vec}(\mathbf{X}_{t-1})\\1\end{pmatrix}\in\mathbb{R}^{(mn+1)\times 1},\quad\mathbf{\Psi}=(\mathbf{\Phi}\vdots\operatorname{vec}(\mathbf{D}))\in\mathbb{R}^{mn\times(mn+1)}.
$$
Here, the matrix $\mathbf{\Phi}$ serves as the comprehensive autoregressive coefficient matrix, while $\operatorname{vec}(\mathbf{D})$ represents the vectorized intercept term. 

The conditional least squares estimate $\hat{\mathbf{\Psi}}_{\mathrm{cls}}$ can be obtained by minimizing the following objective function:
$$
\begin{aligned}
\hat{\mathbf{\Psi}}_{\mathrm{cls}}&=\arg\min_{\mathbf{\Psi}}\sum_{t=2}^T \left\|\operatorname{vec}(\mathbf{X}_t)-\mathbb{E}(\operatorname{vec}(\mathbf{X}_t)|\operatorname{vec}(\mathbf{X}_{t-1}))\right\|_F^2\\
&=\arg\min_{\mathbf{\Psi}}\sum_{t=2}^T \left\|\operatorname{vec}(\mathbf{X}_t)-\mathbf{\Psi}\operatorname{vec}(\mathbf{Z}_{t-1})\right\|_F^2,
\end{aligned}
$$
which yields the closed-form solution
$$
\hat{\mathbf{\Psi}}_{\mathrm{cls}}=\sum_{t=2}^T \operatorname{vec}(\mathbf{X}_t)\operatorname{vec}(\mathbf{Z}_{t-1})^\top \left(\sum_{t=2}^T \operatorname{vec}(\mathbf{Z}_{t-1})\operatorname{vec}(\mathbf{Z}_{t-1})^\top\right)^{-1}.
$$
Once the parameter matrix estimate $\hat{\mathbf{\Psi}}_{\mathrm{cls}}$ is obtained, we can partition it to extract the unconstrained preliminary estimates $\hat{\mathbf{\Phi}}_{\mathrm{cls}}$ and $\hat{\mathbf{D}}_{\mathrm{cls}}$. Recalling the definition $\smash{\mathbf{\Psi} = (\mathbf{\Phi}\vdots\operatorname{vec}(\mathbf{D}))}$, the matrix $\hat{\mathbf{\Phi}}_{\mathrm{cls}}$ is directly extracted as the first $mn$ columns of $\hat{\mathbf{\Psi}}_{\mathrm{cls}}$. Simultaneously, the vectorized intercept estimate $\operatorname{vec}(\hat{\mathbf{D}}_{\mathrm{cls}})$ corresponds to the last (i.e., the $(mn+1)$-th) column of $\hat{\mathbf{\Psi}}_{\mathrm{cls}}$, mathematically denoted as $\smash{\hat{\mathbf{\Psi}}_{\mathrm{cls}, 1:mn,\, mn+1}}$.

Let $\hat{\mathbf{F}}_{\mathrm{cls},k,s}$ represent the $(k,s)$-th sub-block matrix of $\hat{\mathbf{\Phi}}_{\mathrm{cls}}$, that is, 
\[
\hat{\mathbf{F}}_{\mathrm{cls},k,s} = \hat{\mathbf{\Phi}}_{\mathrm{cls}, (k-1)m+1 : km, (s-1)m+1 : sm} \in \mathbb{R}^{m \times m},
\]
referring to the block matrix of $\hat{\mathbf{\Phi}}_{\mathrm{cls}}$ from the $(k-1)m+1$-th row to the $km$-th row and from the $(s-1)m+1$-th column to the $sm$-th column. Consequently, it holds that $\hat{\mathbf{F}}_{\mathrm{cls},k,s} = (\mathbf{e}_{n,k}^\top \otimes \mathbf{I}_{m \times m})\hat{\mathbf{\Phi}}_{\mathrm{cls}}(\mathbf{e}_{n,s} \otimes \mathbf{I}_{m \times m})$, where $1 \leq k, s \leq n$, and $k, s$ are the row and column indices of the sub-block matrix. Accordingly, the explicit PROJs \((\hat{\mathbf{A}}_{\mathrm{PROJ}}, \hat{\mathbf{B}}_{\mathrm{PROJ}}, \hat{\mathbf{C}}_{\mathrm{PROJ}}, \hat{\mathbf{D}}_{\mathrm{PROJ}})\) are formulated as follows:
\begin{align*}
\hat{\mathbf{A}}_{\mathrm{PROJ}} &= \frac{1}{n} \sum_{s=1}^n \left( \hat{\mathbf{F}}_{\mathrm{cls},s,s} - \operatorname{Diag}(\hat{\mathbf{F}}_{\mathrm{cls},s,s}) \right), \\[6pt]
\hat{\mathbf{B}}_{\mathrm{PROJ}} &= (\hat{b}_{k,s})_{n\times n}, \quad \text{where } \hat{b}_{k,s} = 
\begin{cases} 
\displaystyle \frac{1}{m} \mathbf{1}_{m\times1}^\top\operatorname{diag}(\hat{\mathbf{F}}_{\mathrm{cls},k,s}), & k \neq s, \\[4pt]
0, & k = s, 
\end{cases} \\[6pt]
\hat{\mathbf{C}}_{\mathrm{PROJ}} &= \mathcal{R}_{m,n}\big( \operatorname{diag}(\hat{\mathbf{\Phi}}_{\mathrm{cls}}) \big), \\[6pt]
\hat{\mathbf{D}}_{\mathrm{PROJ}} &= \mathcal{R}_{m,n}\big( \hat{\mathbf{\Psi}}_{\mathrm{cls}, 1:mn,\, mn+1} \big).
\end{align*}

\subsection{Iterative Conditional Least Squares Estimation} \label{ICLSE}

Inspired by the alternating least squares approach utilized for continuous additive matrix autoregressions \citep{Zhang2024}, we propose an Iterative Conditional Least Squares Estimation (ICLSE) procedure to estimate the parameter matrices \(\mathbf{A}, \mathbf{B}, \mathbf{C}\), and \(\mathbf{D}\) in Model (III). This method relies on an alternating block coordinate descent algorithm adapted for linear equality constraints. The primary objective is to minimize the conditional sum of squared residuals:
\begin{equation}\label{eq:iclse_obj}
\min_{{\mathbf{A}},{\mathbf{B}},{\mathbf{C}},{\mathbf{D}}} \mathcal{Q}(\mathbf{A},\mathbf{B},\mathbf{C},\mathbf{D}) = \frac{1}{T-1}\sum_{t=2}^T\|\mathbf{X}_t- {\mathbf{A}}\mathbf{X}_{t-1}-\mathbf{X}_{t-1} {\mathbf{B}}^\top- {\mathbf{C}}\odot\mathbf{X}_{t-1}- {\mathbf{D}}\|_F^2,
\end{equation}
subject to the parameter identification constraints \(\operatorname{diag}(\mathbf{A}) = \mathbf{0}_{m \times 1}\) and \(\operatorname{diag}(\mathbf{B}) = \mathbf{0}_{n \times 1}\).

To mathematically handle these zero-diagonal restrictions within the optimization framework, we reformulate them as equivalent trace constraints: \(\operatorname{tr}( {\mathbf{A}} \mathbf{e}_{m,k} \mathbf{e}_{m,k}^\top) = 0\) for \(k=1,\ldots,m\) and \(\operatorname{tr}( {\mathbf{B}} \mathbf{e}_{n,s} \mathbf{e}_{n,s}^\top) = 0\) for \(s=1,\ldots,n\). By introducing the respective Lagrange multiplier vectors \(\boldsymbol{\mu} \in \mathbb{R}^m\) and \(\boldsymbol{\nu} \in \mathbb{R}^n\), the constrained optimization is naturally converted into an unconstrained Lagrangian function \(\mathcal{L}(\mathbf{A},\mathbf{B},\mathbf{C},\mathbf{D},\boldsymbol{\mu},\boldsymbol{\nu})\).

Instead of solving the massive global First-Order Conditions (FOC) simultaneously, our algorithm iteratively updates each parameter block by minimizing the Lagrangian with respect to one block while holding the others fixed, akin to the block coordinate descent methodology \citep{Zhang2024}. Initialized with the projection estimates \((\hat{\mathbf{A}}_{0}, \hat{\mathbf{B}}_{0}, \hat{\mathbf{C}}_{0}, \hat{\mathbf{D}}_{0})\), the \(r\)-th iteration efficiently alternates through the following closed-form updates:

\textbf{Step 1: Update for \(\mathbf{A}\) and \(\boldsymbol{\mu}\).} Conditional on the previous estimates, the optimal row effects and their corresponding multipliers are updated by:
\begin{align*}
\hat{\boldsymbol{\mu}}^r &= 2\operatorname{diag}\!\left[ \left\{\frac{1}{T-1}\sum_{t=2}^T(\mathbf{X}_t-\mathbf{X}_{t-1}\hat{\mathbf{B}}_{r-1}^\top-\hat{\mathbf{C}}_{r-1}\odot\mathbf{X}_{t-1}-\hat{\mathbf{D}}_{r-1})\mathbf{X}_{t-1}^\top\right\} \Upsilon_1^{-1} \right] \oslash \operatorname{diag}\!\left(\Upsilon_1^{-1}\right),\\
\hat{\mathbf{A}}_r &= \left\{ \frac{1}{T-1}\sum_{t=2}^T(\mathbf{X}_t-\mathbf{X}_{t-1}\hat{\mathbf{B}}_{r-1}^\top-\hat{\mathbf{C}}_{r-1}\odot\mathbf{X}_{t-1}-\hat{\mathbf{D}}_{r-1})\mathbf{X}_{t-1}^\top - \sum_{k=1}^m\frac{\hat{\boldsymbol{\mu}}_k^{r}}{2}\mathbf{e}_{m,k} \mathbf{e}_{m,k}^\top \right\} \Upsilon_1^{-1},
\end{align*}
where \(\Upsilon_1 = \frac{1}{T-1}\sum_{t=2}^T\mathbf{X}_{t-1}\mathbf{X}_{t-1}^\top\) denotes the uncentered empirical autocovariance matrix of the lagged observations.

\textbf{Step 2: Update for \(\mathbf{B}\) and \(\boldsymbol{\nu}\).} Symmetrically, the column effects are updated by:
\begin{align*}
\hat{\boldsymbol{\nu}}^r &= 2\operatorname{diag}\!\left[ \Upsilon_2^{-1} \left\{ \frac{1}{T-1}\sum_{t=2}^T\mathbf{X}_{t-1}^\top(\mathbf{X}_t-\hat{\mathbf{A}}_r\mathbf{X}_{t-1}-\hat{\mathbf{C}}_{r-1}\odot\mathbf{X}_{t-1}-\hat{\mathbf{D}}_{r-1}) \right\} \right] \oslash \operatorname{diag}\!\left(\Upsilon_2^{-1}\right),\\
\hat{\mathbf{B}}_r^\top &= \Upsilon_2^{-1} \left\{ \frac{1}{T-1}\sum_{t=2}^T\mathbf{X}_{t-1}^\top(\mathbf{X}_t-\hat{\mathbf{A}}_r\mathbf{X}_{t-1}-\hat{\mathbf{C}}_{r-1}\odot\mathbf{X}_{t-1}-\hat{\mathbf{D}}_{r-1}) - \sum_{s=1}^n\frac{\hat{\boldsymbol{\nu}}_s^{r}}{2}\mathbf{e}_{n,s}\mathbf{e}_{n,s}^\top \right\},
\end{align*}
where \(\Upsilon_2 = \frac{1}{T-1}\sum_{t=2}^T\mathbf{X}_{t-1}^\top\mathbf{X}_{t-1}\).

\textbf{Step 3: Update for \(\mathbf{C}\) and \(\mathbf{D}\).} The unconstrained idiosyncratic lags and intercepts are updated independently via:
\begin{align*}
\hat{\mathbf{C}}_r &= \left\{ \sum_{t=2}^T (\mathbf{X}_t - \hat{\mathbf{A}}_r\mathbf{X}_{t-1} - \mathbf{X}_{t-1}\hat{\mathbf{B}}_r^\top - \hat{\mathbf{D}}_{r-1}) \odot \mathbf{X}_{t-1} \right\} \oslash \left( \sum_{t=2}^T \mathbf{X}_{t-1} \odot \mathbf{X}_{t-1} \right),\\[1ex]
\hat{\mathbf{D}}_r &= \frac{2}{T-1}\sum_{t=2}^T (\mathbf{X}_t - \hat{\mathbf{A}}_r\mathbf{X}_{t-1} - \mathbf{X}_{t-1}\hat{\mathbf{B}}_r^\top - \hat{\mathbf{C}}_r\odot\mathbf{X}_{t-1}).
\end{align*}

This alternating procedure is systematically repeated until convergence. In practical computation, the ICLSE algorithm typically achieves rapid convergence within a few iterations. Regarding the stopping criterion, we adopt the convergence thresholds utilized in recent matrix-variate literature \cite{Xu2025}. Specifically, the iteration is terminated when the squared Frobenius norm of the changes in all estimated matrices falls below a pre-specified tolerance level. The explicit stopping conditions are defined as:
\[
\begin{aligned}
\|\hat{\mathbf{A}}_{r} - \hat{\mathbf{A}}_{r-1}\|_{F}^{2} &< 1 \times 10^{-9}, \quad & 
\|\hat{\mathbf{B}}_{r} - \hat{\mathbf{B}}_{r-1}\|_{F}^{2} &< 1 \times 10^{-9}, \\[1ex]
\|\hat{\mathbf{C}}_{r} - \hat{\mathbf{C}}_{r-1}\|_{F}^{2} &< 1 \times 10^{-9}, \quad & 
\|\hat{\mathbf{D}}_{r} - \hat{\mathbf{D}}_{r-1}\|_{F}^{2} &< 1 \times 10^{-9}.
\end{aligned}
\]
Once all the above inequalities are simultaneously satisfied, the optimization algorithm halts, yielding the final parameter estimates \((\hat{\mathbf{A}}, \hat{\mathbf{B}}, \hat{\mathbf{C}}, \hat{\mathbf{D}})\).

\section{Theory}
In this section, we investigate the asymptotic distributions of the projection estimator and the ICLSE for Model (III). To establish the asymptotic theory, we first define the probabilistic properties of the innovation sequence for the proposed model under a generalized framework allowing for contemporaneous cross-sectional dependence. Recall that \(\mathbf{\Phi} = \mathbf{I}_{n \times n} \otimes \mathbf{A} + \mathbf{B} \otimes \mathbf{I}_{m \times m} + \operatorname{diag}(\operatorname{vec}(\mathbf{C}))\) denotes the comprehensive parameter matrix.

\begin{theorem}\label{t:4.1}
Let \(\{\mathbf{X}_t\}_{t \in \mathbb{Z}}\) be a matrix sequence generated by an \(m \times n\)-dimensional Add-MINAR(1) model with coefficient matrices \(\mathbf{A} \in \mathbb{R}_{\geq 0}^{m \times m}\), \(\mathbf{B} \in \mathbb{R}_{\geq 0}^{n \times n}\), and \(\mathbf{C} \in \mathbb{R}_{\geq 0}^{m \times n}\) satisfying \(\rho(\mathbf{\Phi}) < 1\). Assume the contemporaneous innovation \(\mathbf{E}_t\) possesses an expectation \(\mathbb{E}(\mathbf{E}_t) = \mathbf{D} \in \mathbb{R}_{\geq 0}^{m \times n} \setminus \{\mathbf{0}_{m \times n}\}\) and a general positive definite covariance matrix \(\Sigma_{\mathbf{E}} \coloneqq \operatorname{Cov}(\operatorname{vec}(\mathbf{E}_t), \operatorname{vec}(\mathbf{E}_t))\). Then for \(t \in \mathbb{Z}\),
\[
\mathbf{U}_t \coloneqq \mathbf{X}_t - \mathbf{A} \mathbf{X}_{t-1} - \mathbf{X}_{t-1} \mathbf{B}^\top - \mathbf{C} \odot \mathbf{X}_{t-1} - \mathbf{D}
\]
defines a matrix white noise sequence; that is, \(\{\mathbf{U}_t\}_{t \in \mathbb{Z}}\) is stationary, and for all \(t \in \mathbb{Z}\), \(\mathbb{E}(\mathbf{U}_t) = \mathbf{0}_{m \times n}\) and
\[
\mathbb{E}\bigl(\operatorname{vec}(\mathbf{U}_{t_1}) \operatorname{vec}(\mathbf{U}_{t_2})^\top\bigr) = 
\begin{cases} 
\operatorname{diag}\bigl\{\mathbf{\Phi}(\mathbf{I}_{mn \times mn} - \mathbf{\Phi})^{-1} \operatorname{vec}(\mathbf{D})\bigr\} + \Sigma_{\mathbf{E}}, & t_1 = t_2; \\
\mathbf{0}_{mn \times mn}, & t_1 \neq t_2.
\end{cases}
\]
Moreover, for any \(t_1 < t\), \(\mathbf{U}_t\) is uncorrelated with \(\mathbf{X}_{t_1}\).
\end{theorem}

The proof of Theorem \ref{t:4.1} is provided in the Appendix. The following corollary follows from Theorem \ref{t:4.1}.

\begin{corollary}\label{corollary:4.1}
If the elements of the innovation matrix \(\mathbf{E}_t\) are strictly cross-sectionally independent (e.g., consisting of independent Poisson variables as originally defined), the covariance matrix simplifies to \(\Sigma_{\mathbf{E}} = \operatorname{diag}(\operatorname{vec}(\mathbf{D}))\). Under this strict independence condition, the covariance matrix of \(\mathbf{U}_t\) perfectly reduces to a purely diagonal structure:
\[
\mathbb{E}\bigl(\operatorname{vec}(\mathbf{U}_{t_1}) \operatorname{vec}(\mathbf{U}_{t_2})^\top\bigr) = 
\begin{cases} 
\operatorname{diag}\bigl\{(\mathbf{I}_{mn \times mn} - \mathbf{\Phi})^{-1} \operatorname{vec}(\mathbf{D})\bigr\}, & t_1 = t_2; \\
\mathbf{0}_{mn \times mn}, & t_1 \neq t_2.
\end{cases}
\]
\end{corollary}

Treating Model \eqref{equ:Add-MINAR to MGINAR} as a standard MGINAR(1) process, we establish the asymptotic normality of the conditional least squares estimators by invoking the fundamental regularity conditions outlined by \citet{Latour1997}. Specifically, the following standard conditions are required:
\begin{itemize}
    \item[\bf (I)] \textbf{Stationarity and Ergodicity:} The spectral radius of the comprehensive autoregressive coefficient matrix must satisfy $\rho(\mathbf{\Phi}) < 1$. This guarantees that the process $\{\mathbf{X}_t\}_{t \in \mathbb{Z}}$ admits a unique, strictly stationary, and ergodic solution.
    \item[\bf (II)] \textbf{Innovation Regularity:} The vectorized contemporaneous innovation sequence, $\operatorname{vec}(\mathbf{E}_t)$, comprises independent and identically distributed random vectors with a finite mean $\operatorname{vec}(\mathbf{D})$ and a finite, positive-definite covariance matrix $\Sigma_{\mathbf{E}}$. Additionally, $\mathbf{E}_t$ is independent of past observations $\mathbf{X}_s$ for all $s < t$.
    \item[\bf (III)] \textbf{Non-singularity:} The unconditional covariance matrix of the observation sequence, defined as $\Gamma_0 \coloneqq \operatorname{Cov}(\operatorname{vec}(\mathbf{X}_t), \operatorname{vec}(\mathbf{X}_t))$, is non-singular, ensuring the existence of its inverse $\Gamma_0^{-1}$. This property is intrinsically guaranteed by the positive-definiteness of $\Sigma_{\mathbf{E}}$ specified in Condition (II)\citep{Latour1997}.
\end{itemize}

Provided these fundamental conditions hold, the unconstrained global estimator $\hat{\mathbf{\Phi}}_{\mathrm{cls}}$ is strongly consistent and asymptotically normal. Specifically, following the classical derivation framework of \citet{Latour1997}, the asymptotic distribution of $\operatorname{vec}(\hat{\mathbf{\Phi}}_{\mathrm{cls}})$ is rigorously established as:
\begin{equation}\label{eq:asymptotic_phi}
\sqrt{T}\,\operatorname{vec}(\hat{\mathbf{\Phi}}_{\mathrm{cls}} - \mathbf{\Phi})\xrightarrow{d}\mathbb{N}\bigl(\mathbf{0}_{m^2n^2\times1},\ \Xi\bigr),
\end{equation}
where $\Sigma_{\mathbf{U}} \coloneqq \operatorname{Cov}(\operatorname{vec}(\mathbf{U}_t), \operatorname{vec}(\mathbf{U}_t))$, and $\Xi = \Gamma_0^{-1}\otimes\Sigma_{\mathbf{U}}$.

To intuitively illustrate the covariance structure of the matrix white noise sequence \(\mathbf{U}_t\) established  under the above regularity conditions for asymptotic normality (i.e., zero-diagonal constraints, stationarity \(\rho(\mathbf{\Phi}) < 1\), and non-singularity), consider a low-dimensional bivariate setting where \(m=2\) and \(n=2\).

\begin{example}\label{ex:4.1}
 Suppose the observation matrix is \(\mathbf{X}_t \in \mathbb{N}_0^{2 \times 2}\), and we deliberately specify the interactive coefficient matrices \(\mathbf{A}\) and \(\mathbf{B}\) such that they satisfy the zero-diagonal constraints but are strictly singular (non-invertible):
\[
\mathbf{A} = \begin{bmatrix} 0 & 0.5 \\ 0 & 0 \end{bmatrix}, \quad \mathbf{B} = \begin{bmatrix} 0 & 0 \\ 0.5 & 0 \end{bmatrix}, \quad \mathbf{C} = \begin{bmatrix} 0.5 & 0.5 \\ 0.5 & 0.5 \end{bmatrix}.
\]
Note that \(\det(\mathbf{A}) = 0\) and \(\det(\mathbf{B}) = 0\), rendering both spatial/cross-sectional matrices non-invertible. Assume the Poisson innovation parameter matrix is uniformly given by:
\[
\mathbf{D} = \begin{bmatrix} 1 & 1 \\ 1 & 1 \end{bmatrix} \implies \operatorname{vec}(\mathbf{D}) = (1, 1, 1, 1)^\top.
\]
The comprehensive autoregressive coefficient matrix \(\mathbf{\Phi}\) is constructed as:
\[
\mathbf{\Phi} = \mathbf{I}_{2 \times 2} \otimes \mathbf{A} + \mathbf{B} \otimes \mathbf{I}_{2 \times 2} + \operatorname{diag}(\operatorname{vec}(\mathbf{C})) = \begin{bmatrix} 0.5 & 0.5 & 0 & 0 \\ 0 & 0.5 & 0 & 0 \\ 0.5 & 0 & 0.5 & 0.5 \\ 0 & 0.5 & 0 & 0.5 \end{bmatrix}.
\]
It is straightforward to verify that the spectral radius is \(\rho(\mathbf{\Phi}) = 0.5 < 1\), thereby strictly satisfying the stationarity condition required for asymptotic normality. According to the formula in Corollary \ref{corollary:4.1}, we subsequently compute the multiplier matrix:
\[
(\mathbf{I}_{4 \times 4} - \mathbf{\Phi})^{-1} = \begin{bmatrix} 0.5 & -0.5 & 0 & 0 \\ 0 & 0.5 & 0 & 0 \\ -0.5 & 0 & 0.5 & -0.5 \\ 0 & -0.5 & 0 & 0.5 \end{bmatrix}^{-1} = \begin{bmatrix} 2 & 2 & 0 & 0 \\ 0 & 2 & 0 & 0 \\ 2 & 4 & 2 & 2 \\ 0 & 2 & 0 & 2 \end{bmatrix}.
\]
Thus, when \(t_1 = t_2\), the unconditional (contemporaneous) covariance matrix of the vectorized white noise evaluates to:
\[
\begin{aligned}
\mathbb{E}\bigl(\operatorname{vec}(\mathbf{U}_t) \operatorname{vec}(\mathbf{U}_t)^\top\bigr) &= \operatorname{diag}\bigl\{ (\mathbf{I}_{4 \times 4} - \mathbf{\Phi})^{-1} \operatorname{vec}(\mathbf{D}) \bigr\} \\
&= \operatorname{diag}\left\{ \begin{bmatrix} 2 & 2 & 0 & 0 \\ 0 & 2 & 0 & 0 \\ 2 & 4 & 2 & 2 \\ 0 & 2 & 0 & 2 \end{bmatrix} \begin{pmatrix} 1 \\ 1 \\ 1 \\ 1 \end{pmatrix} \right\} \\
&= \operatorname{diag}(4, 2, 10, 4).
\end{aligned}
\]
Furthermore, to characterize the stationary variability of the observed matrix sequence, we compute the covariance matrix \(\Gamma_0 := \operatorname{Cov}(\operatorname{vec}(\mathbf{X}_t), \operatorname{vec}(\mathbf{X}_t))\), which satisfies the discrete algebraic Lyapunov equation \(\Gamma_0 = \mathbf{\Phi} \Gamma_0 \mathbf{\Phi}^\top + \Sigma_{\mathbf{U}}\). Solving this via the vectorization property \(\operatorname{vec}(\Gamma_0) = (\mathbf{I}_{4 \times 4} - \mathbf{\Phi} \otimes \mathbf{\Phi})^{-1} \operatorname{vec}(\Sigma_{\mathbf{U}})\), we obtain the stationary covariance:
\[
\Gamma_0 \approx \begin{bmatrix} 
6.81 & 0.89 & 3.56 & 1.48 \\ 
0.89 & 2.67 & 0.59 & 0.89 \\ 
3.56 & 0.59 & 23.60 & 3.56 \\ 
1.48 & 0.89 & 3.56 & 6.81 
\end{bmatrix}.
\]
\end{example}

The above example illustrates an important property: despite both row and column interactive matrices \(\mathbf{A}\) and \(\mathbf{B}\) being singular, the contemporaneous covariance matrix of the white noise is \(\operatorname{diag}(4, 2, 10, 4)\) and $\Gamma_0$ is positive definite. This example visually corroborates that the proposed model readily satisfies the non-singularity assumption essential for asymptotic normality, regardless of the invertibility of the interactive coefficient matrices.

\begin{proposition}\label{prop:positive_definite}
Under the stationarity condition \(\rho(\mathbf{\Phi}) < 1\), if the covariance matrix of the innovation sequence \(\Sigma_{\mathbf{E}}\) is positive definite, then the contemporaneous covariance matrix of the matrix white noise \(\Sigma_{\mathbf{U}}\) and the unconditional covariance matrix of the observation sequence \(\Gamma_0\) are both strictly positive definite.
\end{proposition}
The detailed proof of Proposition \ref{prop:positive_definite} is provided in the Appendix.

Since the diagonal elements of the matrix estimators for \(\mathbf{A}\) and \(\mathbf{B}\) are always zero, their asymptotic variances vanish. To avoid redundant content and for ease of exposition, for \(\mathbf{A} \in \mathbb{R}^{m \times m}\), let \(\operatorname{vec}_m(\mathbf{A}) \in \mathbb{R}^{m^2 - m}\) denote the subvector of \(\operatorname{vec}(\mathbf{A})\) excluding its diagonal elements.

\begin{theorem}\label{t:4.2}
For Model (III), suppose that \(\Sigma_{\mathbf{E}}\) is positive definite and finite. Then, as \(T \to \infty\),
\[
\sqrt{T}\,\operatorname{vec}_m(\hat{\mathbf{A}}_{\mathrm{PROJ}} - \mathbf{A}) \xrightarrow{d} \mathbb{N}\bigl(\mathbf{0}_{(m^2-m) \times 1},\; \Gamma_1^* \Xi (\Gamma_1^*)^\top\bigr),
\]
\[
\sqrt{T}\,(\hat{b}_{k,s} - b_{k,s}) \xrightarrow{d} \mathbb{N}\bigl(0,\; \Gamma_2 \Xi \Gamma_2^\top\bigr), \quad (k \neq s),
\]
\[
\sqrt{T}\,\operatorname{vec}(\hat{\mathbf{C}}_{\mathrm{PROJ}} - \mathbf{C}) \xrightarrow{d} \mathbb{N}\bigl(\mathbf{0}_{mn \times 1},\; \Gamma_3 \Xi \Gamma_3^\top\bigr),
\]
\[
\sqrt{T}\,\operatorname{vec}(\hat{\mathbf{D}}_{\mathrm{PROJ}} - \mathbf{D}) \xrightarrow{d} \mathbb{N}\bigl(\mathbf{0}_{mn \times 1},\; \bigl(1+\operatorname{vec}(\boldsymbol{\mu}_{\mathbf{X}})^\top \Gamma_0^{-1}\operatorname{vec}(\boldsymbol{\mu}_{\mathbf{X}})\bigr) \Sigma_{\mathbf{U}}\bigr),
\]
where \(\boldsymbol{\mu}_{\mathbf{X}} \coloneqq \mathbb{E}(\mathbf{X}_t)\), and
\[
\begin{aligned}
\Gamma_{1} &= \left\{\mathbf{I}_{m^2\times m^2} - \sum_{k=1}^m \left\{(\mathbf{e}_{m,k}\otimes \mathbf{e}_{m,k})(\mathbf{e}_{m,k}^\top\otimes \mathbf{e}_{m,k}^\top)\right\}\right\} \nonumber 
\cdot \frac{1}{n} \sum_{s=1}^n \left\{(\mathbf{e}_{n,s}^\top\otimes \mathbf{I}_{m\times m})\otimes(\mathbf{e}_{n,s}^\top\otimes \mathbf{I}_{m\times m})\right\},\\[6pt]
\Gamma_{2} &= \frac{1}{m} \mathbf{1}_{m\times1}^\top \left\{(\mathbf{e}_{m,1}\otimes \mathbf{e}_{m,1},\;\ldots,\; \mathbf{e}_{m,m}\otimes \mathbf{e}_{m,m})^\top\right\} \left\{(\mathbf{e}_{n,s}^\top\otimes \mathbf{I}_{m\times m})\otimes(\mathbf{e}_{n,k}^\top\otimes \mathbf{I}_{m\times m})\right\},\\[6pt]
\Gamma_{3} &= (\mathbf{e}_{mn,1}\otimes\mathbf{e}_{mn,1},\;\ldots,\; \mathbf{e}_{mn,mn}\otimes\mathbf{e}_{mn,mn})^\top,
\end{aligned}
\]
$\Gamma_1^*$ represents the reduced matrix derived by eliminating the rows of $\Gamma_1$ that correspond to the zero-diagonal constraints of $\mathbf{A}$. Specifically, this operation removes the $((k-1)m+k)$-th rows from $\Gamma_1$ for all $k=1, 2, \dots, m$.
\end{theorem}
The proof of Theorem \ref{t:4.2} is given in the Appendix.

Next, we discuss the joint asymptotic distribution of the ICLSE estimator. Let \(\mathbf{\Pi}_1 = \mathbb{E}(\mathbf{X}_t \mathbf{X}_t^\top)\), \(\mathbf{\Pi}_2 = \mathbb{E}(\mathbf{X}_t^\top \mathbf{X}_t)\), and define the matrices  
\(\mathbf{M}_1 \in \mathbb{R}^{m^2 \times m^2}\), \(\mathbf{M}_2 \in \mathbb{R}^{n^2 \times n^2}\), \(\mathbf{M} \in \mathbb{R}^{\ell \times \ell}\) and \(\mathbf{W}_t \in \mathbb{R}^{\ell \times mn}\) \((\ell = (m+n)^2)\) as follows:
\[
\mathbf{M}_1 = \sum_{k=1}^m \frac{(\mathbf{e}_{m,k} \mathbf{e}_{m,k}^\top \mathbf{\Pi}_1^{-1}) \otimes (\mathbf{e}_{m,k} \mathbf{e}_{m,k}^\top)}{\mathbf{e}_{m,k}^\top \mathbf{\Pi}_1^{-1} \mathbf{e}_{m,k}}, \quad
\mathbf{M}_2 = \sum_{s=1}^n \frac{(\mathbf{e}_{n,s} \mathbf{e}_{n,s}^\top) \otimes (\mathbf{e}_{n,s} \mathbf{e}_{n,s}^\top \mathbf{\Pi}_2^{-1})}{\mathbf{e}_{n,s}^\top \mathbf{\Pi}_2^{-1} \mathbf{e}_{n,s}},
\]
\[
\begin{aligned}
~~\mathbf{M} 
= \operatorname{Blk-diag}\bigl(\mathbf{M}_1 - \mathbf{I}_{m^2 \times m^2},\; \mathbf{M}_2 - \mathbf{I}_{n^2 \times n^2},\; \mathbf{I}_{mn \times mn},\; \mathbf{I}_{mn \times mn}\bigr), 
\quad\mathbf{W}_t 
= \begin{bmatrix}
\mathbf{X}_t \otimes \mathbf{I}_{m \times m} \\[2pt]
\mathbf{I}_{n \times n} \otimes \mathbf{X}_t^\top \\[2pt]
\operatorname{diag}(\operatorname{vec}(\mathbf{X}_t)) \\[2pt]
\mathbf{I}_{n \times n} \otimes \mathbf{I}_{m \times m}
\end{bmatrix}.
\end{aligned}
\]
Here, \(\operatorname{Blk-diag}\) denotes the operation that converts the enclosed elements into a block diagonal matrix. For \(\mathbf{B} \in \mathbb{R}^{n \times n}\), let \(\operatorname{vec}_n(\mathbf{B}) \in \mathbb{R}^{n^2 - n}\) denote the vectorization of \(\mathbf{B}\) with its diagonal elements removed. For any matrix \(\mathbf{V} \in \mathbb{R}^{\ell \times \ell}\), let \(\mathbf{V}^* \in \mathbb{R}^{(\ell - m - n) \times \ell}\) denote the reduced matrix obtainedby removing the \((m(k-1)+k)\)-th and \((m^2 + n(s-1) + s)\)-th rows from \(\mathbf{V}\), where \(k = 1, 2, \dots, m\) and \(s = 1, 2, \dots, n\).

\begin{theorem}\label{t:4.3}
For Model (III), suppose that \(\Sigma_{\mathbf{E}}\) is positive definite and finite. Then, as \(T \to \infty\),
\[
\sqrt{T}
\begin{bmatrix}
\operatorname{vec}_m(\hat{\mathbf{A}} - \mathbf{A}) \\
\operatorname{vec}_n(\hat{\mathbf{B}}^\top - \mathbf{B}^\top) \\
\operatorname{vec}(\hat{\mathbf{C}} - \mathbf{C}) \\
\operatorname{vec}(\hat{\mathbf{D}} - \mathbf{D})
\end{bmatrix}
\xrightarrow{d} \mathbb{N}\left(\mathbf{0}_{(\ell - m - n) \times 1}, \mathcal{X}\right),
\]
where
\[
\mathcal{X} = \left\{(\mathbf{H}^\top)^* ((\mathbf{H}^\top)^*)^\top\right\}^{-1}
\left\{(\mathbf{H}^\top)^* \mathbf{M} \mathbb{E}(\mathbf{W}_t \Sigma_{\mathbf{U}} \mathbf{W}_t^\top) \mathbf{M}^\top ((\mathbf{H}^\top)^*)^\top\right\}
\left\{(\mathbf{H}^\top)^* ((\mathbf{H}^\top)^*)^\top\right\}^{-1},
\]
\[
\begin{aligned}
~\mathbf{H} = \mathbf{M} \mathbb{E}(\mathbf{W}_t \mathbf{W}_t^\top) - \operatorname{Blk-diag}\left(\mathbf{M}_1(\mathbf{\Pi}_1 \otimes \mathbf{I}_{m \times m}),\; \mathbf{M}_2(\mathbf{I}_{n \times n} \otimes \mathbf{\Pi}_2),\; \mathbf{0}_{mn \times mn},\; \mathbf{0}_{mn \times mn}\right).
\end{aligned}
\]
\end{theorem}
The proof of Theorem \ref{t:4.3} is given in the Appendix.

\section{Experiments}
\subsection{Simulation}
In this section, we investigate the estimation performance of the PROJ and ICLSE for the Add-MINAR(1) model. The estimation efficacy of the two parameter estimation methods is demonstrated under different dimensions and sample sizes.

For the true parameter matrices \(\mathbf{A}\), \(\mathbf{B}\), and \(\mathbf{C}\) in Model (III), the elements are generated to be greater than 0 and less than 1, while satisfying \(\rho(\mathbf{I}_{n \times n} \otimes \mathbf{A} + \mathbf{B} \otimes \mathbf{I}_{m \times m} + \operatorname{diag}(\operatorname{vec}(\mathbf{C}))) < 1\), thereby ensuring that the subsequent data generation process is causally stationary. Data \(\mathbf{X}_t\) are then generated according to Model (III). Furthermore, to ensure the uniqueness of the estimated model parameters, the parameter matrices are generated subject to the constraints \(\operatorname{diag}(\mathbf{A}) = \mathbf{0}_{m \times 1}\) and \(\operatorname{diag}(\mathbf{B}) = \mathbf{0}_{n \times 1}\). Regarding the specification of the innovation term \(\mathbf{E}_t\), in order to simulate scenarios where the innovation term follows different levels of complexity, we follow \citet{Cui2026} and consider the following three settings:
\begin{itemize}
  \item[$\bullet$] \textbf{Setting I:} The covariance matrix $\operatorname{Cov}(\operatorname{vec}(\mathbf{E}_t))=\Sigma $ is set to $\Sigma = \mathbf{I}$.
  \item[$\bullet$] \textbf{Setting II:} The covariance matrix $\operatorname{Cov}(\operatorname{vec}(\mathbf{E}_t))=\Sigma$ is diagonal, with its elements following a uniform distribution on (0,1).
  \item[$\bullet$] \textbf{Setting III:} The covariance matrix $\operatorname{Cov}(\operatorname{vec}(\mathbf{E}_t))=\Sigma=\Sigma_c \otimes \Sigma_r$ is randomly generated according to $\Sigma_c=\mathbf{Q}\Lambda \mathbf{Q}^{\top}$, where the eigenvalues in the diagonal matrix $\Lambda$ are the absolute values of i.i.d. standard normal random variates, and the eigenvector matrix $\mathbf{Q}$ is a random orthonormal matrix. $\Sigma_r$ is formed following an identical procedure.
\end{itemize}
 
For the two estimation methods of the Add-MINAR(1) model, we first compare their respective estimation errors under different dimensions and different time series lengths. Here, \(\mathrm{err}_\mathbf{A}\), \(\mathrm{err}_\mathbf{B}\), and \(\mathrm{err}_\mathbf{C}\) denote the errors between the estimated and true coefficient matrices under the projection estimation method, while \(\mathrm{Err}_\mathbf{A}\), \(\mathrm{Err}_\mathbf{B}\), and \(\mathrm{Err}_\mathbf{C}\) denote the corresponding errors under the iterative conditional least squares estimation method. The formulas for calculating each error are given as follows:
\begin{equation}\label{equ:error_exp1}
\begin{aligned}
&\mathrm{err}_{\mathbf{A}} = \log(\|\mathbf{A} - \hat{\mathbf{A}}_{\mathrm{PROJ}}\|_F^{2}),
\ \mathrm{err}_{\mathbf{B}} = \log(\|\mathbf{B} - \hat{\mathbf{B}}_{\mathrm{PROJ}}\|_F^{2}), 
\ \mathrm{err}_{\mathbf{C}} = \log(\|\mathbf{C} - \hat{\mathbf{C}}_{\mathrm{PROJ}}\|_F^{2}), \\
&\mathrm{Err}_{\mathbf{A}} = \log(\|\mathbf{A} - \hat{\mathbf{A}}\|_F^{2}), 
\ \mathrm{Err}_{\mathbf{B}} = \log(\|\mathbf{B} - \hat{\mathbf{B}}\|_F^{2}), 
\ \mathrm{Err}_{\mathbf{C}} = \log(\|\mathbf{C} - \hat{\mathbf{C}}\|_F^{2}).
\end{aligned}
\end{equation}

\begin{figure}[H]
  \centering
  \includegraphics[width=1\textwidth, trim=20 20 20 20, clip]{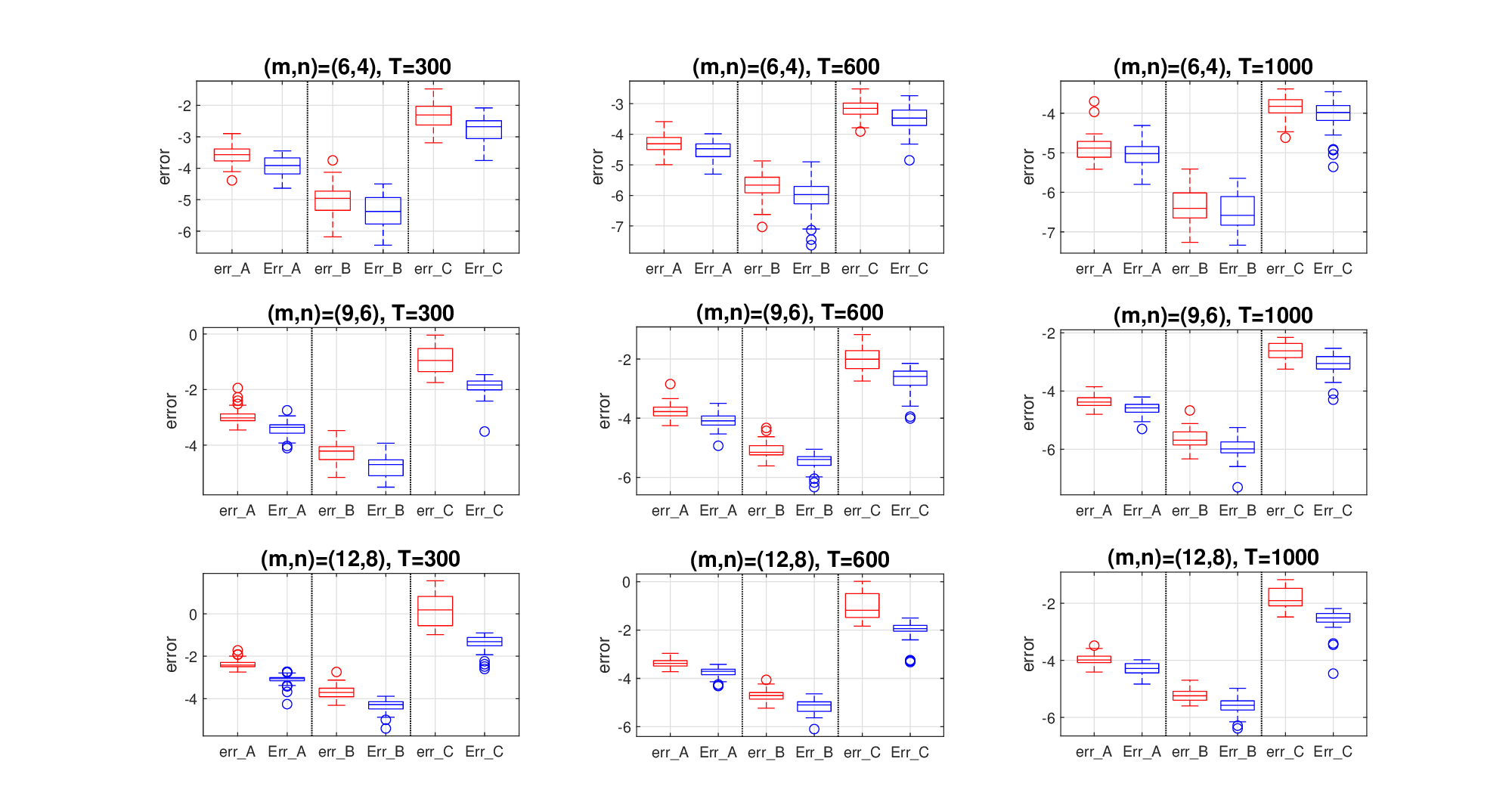}
  \caption{Setting I: Comparison of estimation errors for the coefficient matrices of the Add-MINAR(1) model under the two estimation methods. Red: PROJ results; Blue: ICLSE results.} 
  \label{fig:Add-MINAR-set1_exp1}
\end{figure}

\begin{figure}[htbp]
  \centering
  \includegraphics[width=1\textwidth, trim=20 20 20 20, clip]{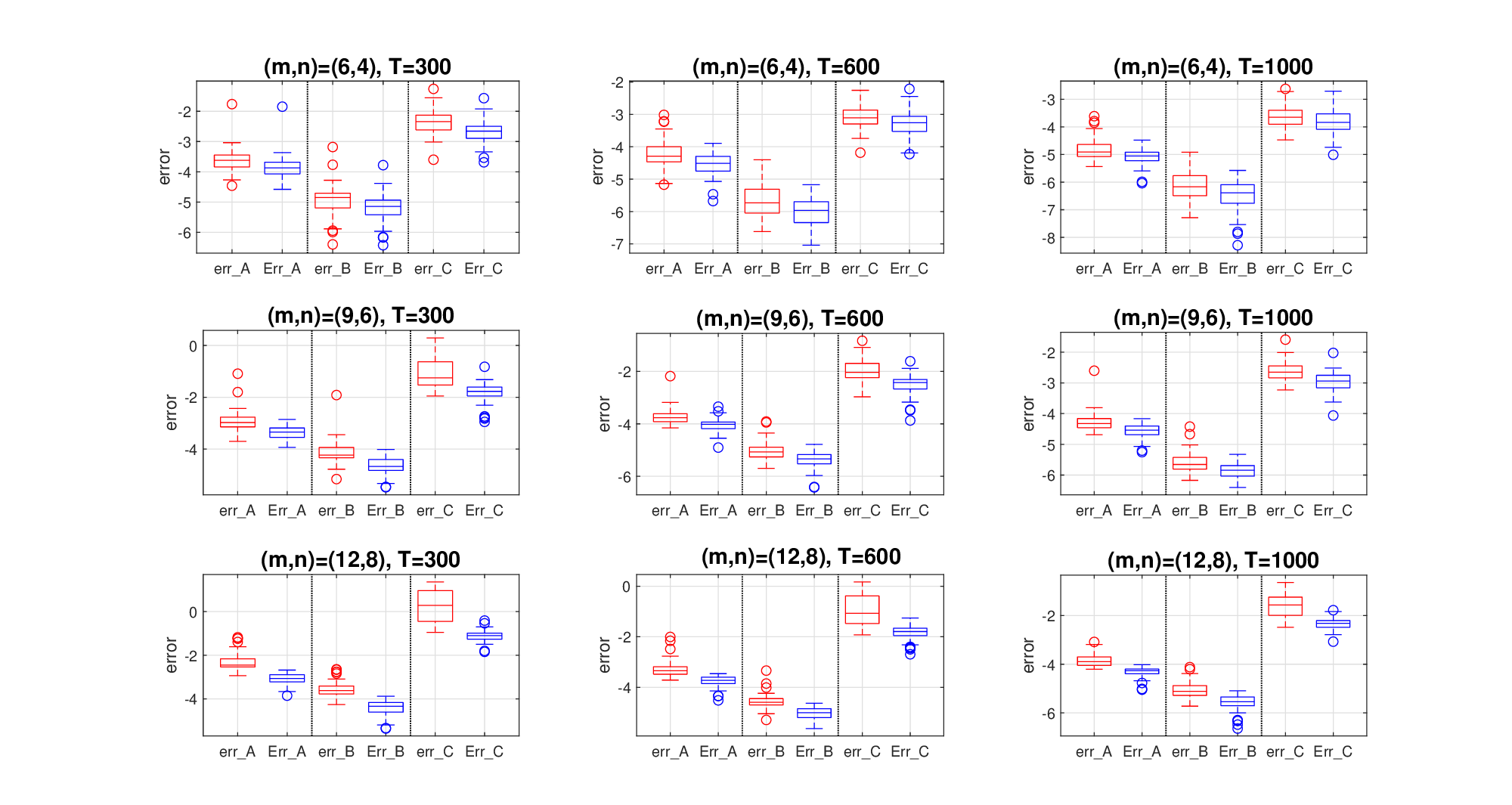}
  \caption{Setting II: Comparison of estimation errors for the coefficient matrices of the Add-MINAR(1) model under the two estimation methods. Red: PROJ results; Blue: ICLSE results.} 
  \label{fig:Add-MINAR-set2_exp1}
\end{figure}

\begin{figure}[htbp]
  \centering
  \includegraphics[width=1\textwidth, trim=20 20 20 20, clip]{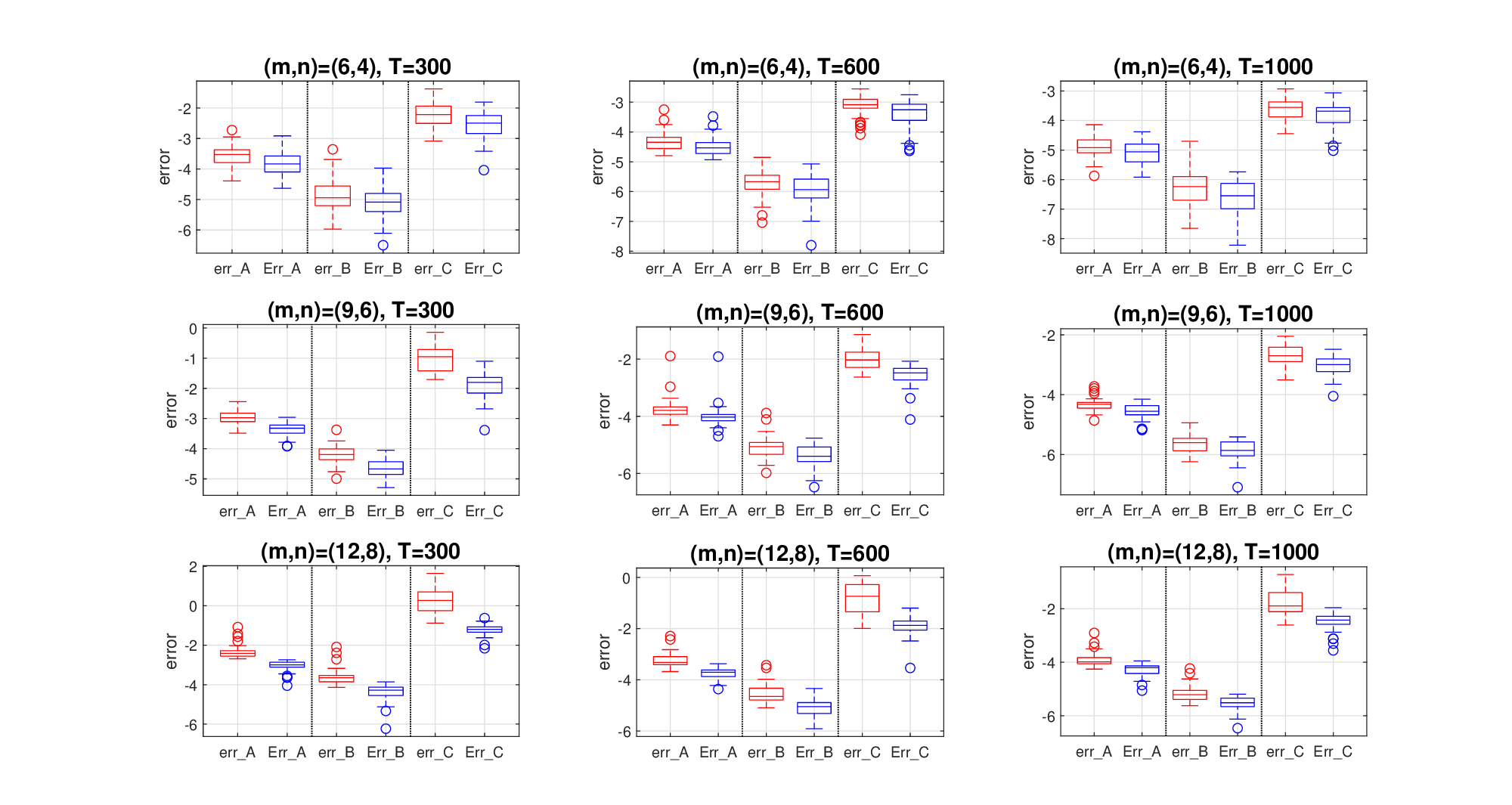}
  \caption{Setting III: Comparison of estimation errors for the coefficient matrices of the Add-MINAR(1) model under the two estimation methods. Red: PROJ results; Blue: ICLSE results.} 
  \label{fig:Add-MINAR-set3_exp1}
\end{figure}

To investigate the performance of the two estimation methods under different settings, we set the dimensions to \((m,n) = (6,4), (9,6), (12,8)\) and the sample sizes to \(T = 300, 600, 1000\). For each combination of dimension and sample size, 100 repeated experiments are conducted. The estimation errors of the coefficient matrices \(\mathbf{A}\), \(\mathbf{B}\), and \(\mathbf{C}\) for both methods are calculated based on \eqref{equ:error_exp1}. The experimental results under different distributions of the parameter matrix \(\mathbf{D}\) are presented in Figures \ref{fig:Add-MINAR-set1_exp1} to \ref{fig:Add-MINAR-set3_exp1}.

Figures \ref{fig:Add-MINAR-set1_exp1} - \ref{fig:Add-MINAR-set3_exp1} display the experimental results under three different settings of the innovation parameter matrix. From the boxplots, it can be observed that under all three experimental settings, the blue boxes (representing ICLSE) are consistently lower than the red boxes (representing PROJ), indicating that ICLSE yields smaller parameter estimation errors and thus higher estimation accuracy. Furthermore, for a fixed data dimension, the estimation errors of both methods decrease significantly as the sample size increases, demonstrating the consistency of the parameter estimates. Under the same sample size, however, the estimation errors of both methods increase with the data dimension. This is because the number of parameters grows substantially as the dimension increases; when the number of parameters approaches or even exceeds the sample size, estimation performance deteriorates and overfitting may occur, which is consistent with the general principles of high-dimensional statistical theory.

\newpage

Furthermore, to further analyze the estimation efficiency of the two methods, we repeatedly computed the estimation errors for both methods under different settings and sample sizes based on 100 replications. The trend of the estimation error as a function of the sample size is presented for each of the three settings and different dimensions. The error is calculated using the following expression:

\[
\mathrm{err}_\Phi = \log\left(\bigl\| \Phi - \bigl( \mathbf{I}_{n\times n} \otimes \hat{\mathbf{A}} + \hat{\mathbf{B}} \otimes \mathbf{I}_{m\times m} + \operatorname{diag}(\operatorname{vec}(\hat{\mathbf{C}})) \bigr) \bigr\|_F^{2} \right).
\]

Under all parameter configurations, a consistent conclusion is observed: for the same dimension, the estimation error of ICLSE is always lower than that of PROJ. As the sample size \(T\) increases, both methods exhibit a clear convergence trend, with the logarithmic error decreasing monotonically, which aligns with the theoretical expectation of consistency. Specifically, in the small-sample regime (\(T < 20,000\)), the estimation errors of both methods decrease rapidly. In the moderate-sample range (\(20,000 < T < 40,000\)), the rate of decline slows, and the curves become flatter with slight fluctuations. As the sample size further increases (\(T > 40,000\)), the error reduction gradually levels off, and in the large-sample regime (\(T > 80,000\)), the errors become essentially stable. In comparison, ICLSE outperforms PROJ for the vast majority of sample sizes, not only exhibiting a lower initial error but also achieving a smaller final convergence value and a more stable convergence trajectory. Therefore, ICLSE demonstrates superior estimation accuracy and numerical stability, making it more suitable for parameter estimation in moderate-to-large sample settings.

\begin{figure}[htbp]
  \centering
  \includegraphics[width=1\textwidth, trim=20 0 20 0, clip]{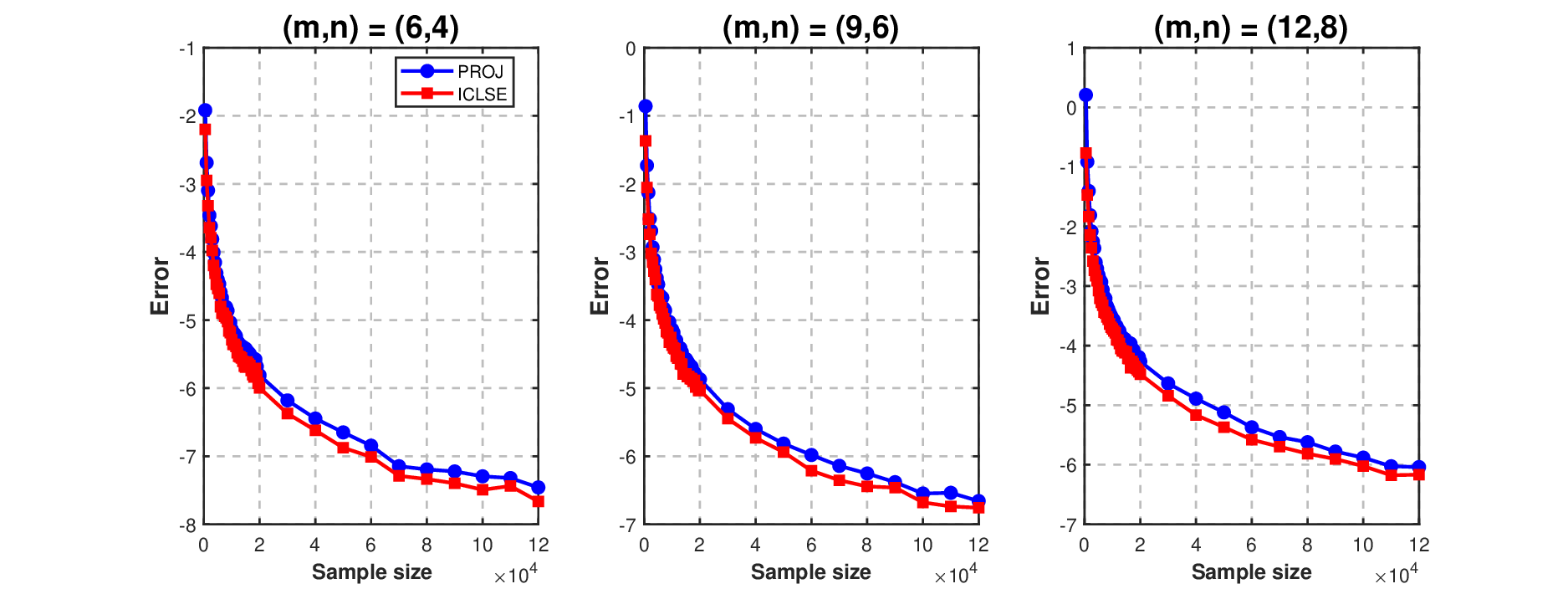}
  \caption{Setting I: Convergence efficiency of the Add-MINAR(1) coefficient matrix \((\mathbf{I}_{n \times n} \otimes \mathbf{A} + \mathbf{B} \otimes \mathbf{I}_{m \times m} + \operatorname{diag}(\operatorname{vec}(\mathbf{C})))\).} 
  \label{fig:Add-MINAR-set1_exp2}
\end{figure}

\subsection{Crime Data Analysis}\label{true data}
In the intersection of crime geography and time series analysis, the interrelationships between different types of crime across different regions have long attracted substantial research attention. Such interrelationships may arise from spatial proximity (criminal activities tend to cluster geographically, forming so-called "crime hot spots"), socioeconomic factors (e.g., poverty, unemployment, low educational attainment, and lack of community resources may render certain types of crime more prevalent in specific areas), crime chain reactions (e.g., associations among prostitution, gambling, and drug-related crimes). To more comprehensively test theories of crime spatial transmission and identify their geographical boundaries, the present study deliberately breaks the assumption of geographical continuity in the spatial dimension and introduces a logic of spatial comparison.

\begin{figure}[H]
\centering
\includegraphics[width=1\textwidth, trim=0 0 0 0, clip]{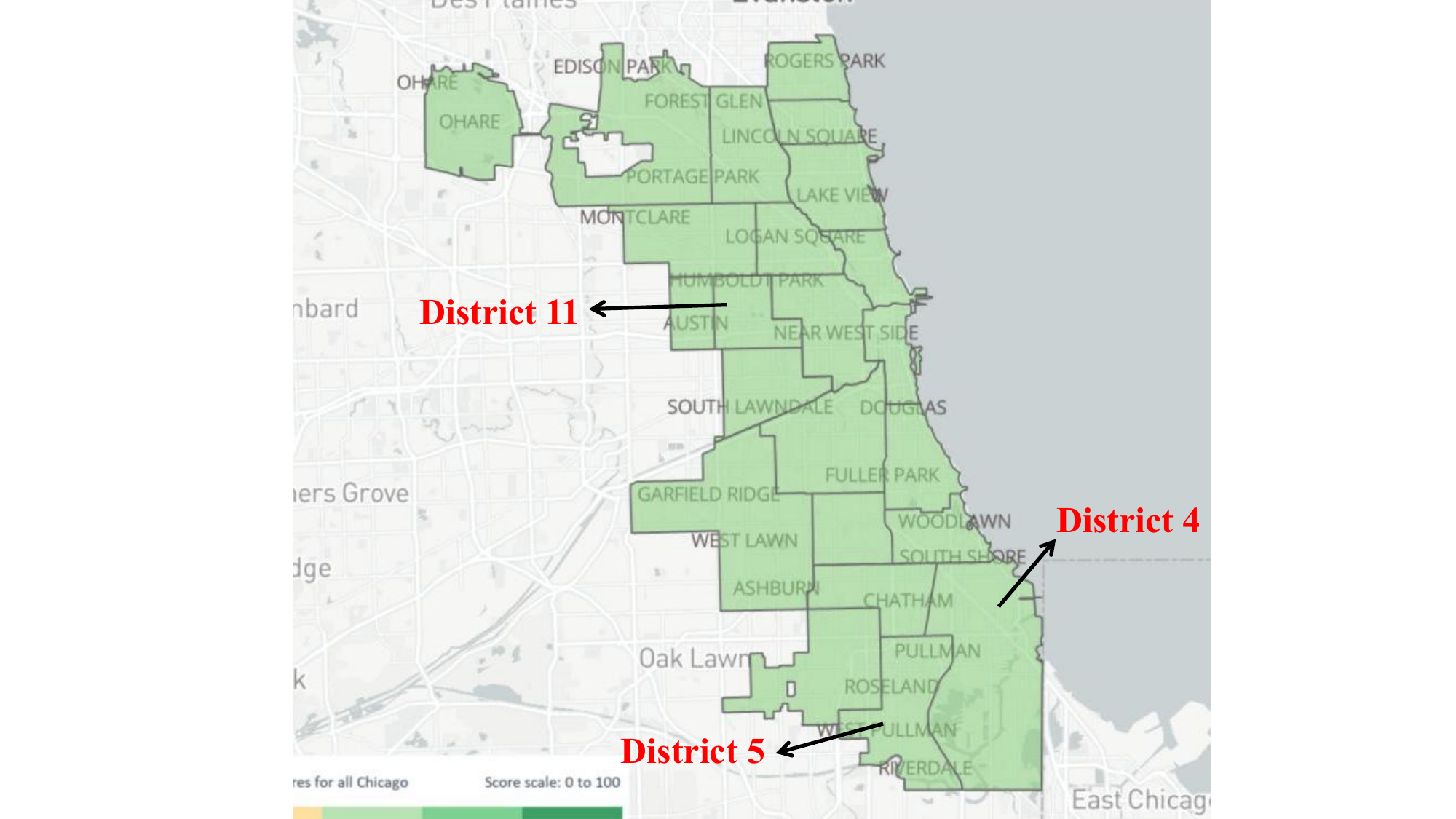}
\caption{Map of police districts in the City of Chicago}
\label{fig:map_Add-MINAR}
\end{figure}

Specifically, three representative police districts in Chicago were selected as the research subjects: Districts 4, 5 and 11. These three districts form a non-contiguous spatial ensemble: Districts 4 and 5 are adjacent to each other and located in the southernmost part of the city, an area characterized by high crime rates; neither of them shares a geographical border with District 11, thereby creating a cross-regional comparative analytical framework (see Figure \ref{fig:map_Add-MINAR} for detailed locations). This design enables the column parameter matrix of the model to simultaneously detect two types of spatial effects: first, possible spatial diffusion or spillover effects of criminal activities between adjacent districts (e.g., between Districts 4 and 5); second, possible crime pattern associations that transcend geographical proximity between non-adjacent districts (e.g., between District 11 and the southern Districts 4 and 5).

Regarding the selection of crime types, this study focuses on three categories characterized by internal heterogeneity and complex social dynamics: Sexual Offenses, Vehicle Theft, and Robbery. Vehicle Theft represents a specific form of general Theft, while Sexual Offenses exhibit entirely different social causes and manifestations. The combination of these three crime types allows the study to move beyond the analysis of associations among homogeneous street crimes and to investigate whether latent concomitant relationships or substitution effects exist among different types of crime.

Based on the above considerations, this study constructs the following data structure using crime data from the City of Chicago. For each observation month \( t \), a \( 3 \times 3 \) crime count matrix \( \mathbf{X}_t \) is defined, where the element \( \mathbf{X}_{t,i,j} \) denotes the number of crimes of the \( i \)-th type occurring in the \( j \)-th police district during month \( t \). This data structure aligns closely with the specifications of the Add-MINAR model: the row parameter matrix \( \mathbf{A} \) can capture mutual influences and concomitant relationships among different crime types; the column parameter matrix \( \mathbf{B} \) can reveal spatial transmission or competition effects of criminal activities across districts; and the parameter matrix \( \mathbf{C} \) can capture the historical dependence characteristics specific to each ``crime type--police district'' combination.

This study employs monthly crime records from 2001.01 to 2023.07, yielding a total of \( 271 \) observations that form the matrix sequence \( \{\mathbf{X}_t\}_{t=1}^{271} \). The data are derived from \url{https://www.payititi.com/opendatasets/show/6962/}. For detailed division rules, please refer to \url{https://www.chicagopolice.org/statistics-data/data-dashboards/sentiment-dashboard/}. The first 251 months of data are used as the training set, and the remaining 20 months as the test set. The advantages of the Add-MINAR model are evaluated by comparing its in-sample and out-of-sample fitting errors against the following five competing models.

\begin{enumerate}[label=(\roman*)]
  \item MGINAR(1): Fits an MGINAR(1) model to $\operatorname{vec}(\mathbf{X}_{t})$.
  \item iINAR(1): Fits nine independent INAR(1) models, each corresponding to a ``crime type $\times$ district'' combination.
  \item MINAR(1)\_ICLSE: Fits an unconstrained rank MINAR(1) model to $\mathbf{X}_{t}$ using ICLSE.
  \item Add-MINAR(1)\_PROJ: Fits an Add-MINAR(1) model to $\mathbf{X}_{t}$ using PROJ.
  \item Add-MINAR(1)\_ICLSE: Fits an Add-MINAR(1) model to $\mathbf{X}_{t}$ using ICLSE.
\end{enumerate}

Table \ref{tab:Add-MINAR_within_sample_fitting_new} presents the in-sample prediction error metrics (\( \mathrm{E}_1 \) to \( \mathrm{E}_4 \)) and the number of parameters for each model. The specific expressions for the error metrics are given in Appendix. In terms of overall prediction accuracy, the Add-MINAR(1)\_ICLSE model performs best, holding a slight advantage in \( \mathrm{E}_1 \) (8192.1987), \( \mathrm{E}_2 \) (11.6915), and \( \mathrm{E}_4 \) (0.2093). It is worth noting that although the Add-MINAR(1)\_PROJ model achieves the lowest \( \mathrm{E}_3 \) (0.8523) among all models, it is markedly inferior to both the MINAR(1) and Add-MINAR(1) models estimated via the ICLSE method in terms of \( \mathrm{E}_1 \), \( \mathrm{E}_2 \), and \( \mathrm{E}_4 \), suggesting that its predictive stability may be insufficient. In contrast, the most structurally complex model, MGINAR(1) (with 90 parameters), exhibits the highest error values across all metrics and thus the poorest predictive performance. The simplest model, iINAR(1) (with 18 parameters), yields smaller errors than MGINAR(1) but still much larger errors than the best-performing models. These findings indicate that, in the context of this study, the Add-MINAR(1) model estimated via the ICLSE method achieves the best balance between prediction accuracy and model complexity, whereas excessive parameterization (as in MGINAR(1)) or overly simplified structure (as in iINAR(1)) both lead to degraded predictive performance.
\begin{table}[htbp]
  \centering
  \caption{Comparison of in-sample fitting performance among five models}
  \label{tab:Add-MINAR_within_sample_fitting_new}
  \begin{tabular}{lccccc}
      \toprule
      Model & E1 & E2 & E3 & E4 & Number of Parameters \\
      \midrule
      MGINAR(1)       & 9787.5930 & 13.7952 & 1.0515 & 0.2474 & 90  \\
      iINAR(1)        & 8766.5851 & 12.4359 & 1.0136 & 0.2243 & \textbf{18} \\
      MINAR(1)\_ICLSE & 8198.6566 & 11.6915 & 0.9461 & 0.2096 & 27  \\
      Add-MINAR(1)\_PROJ & 8892.0183 & 12.5337 & \textbf{0.8523} & 0.2321 & 36  \\
      Add-MINAR(1)\_ICLSE & \textbf{8192.1987} & \textbf{11.6728} & 0.9158 & \textbf{0.2093} & 36  \\
      \bottomrule
  \end{tabular}
\end{table}

\begin{table}[H]
\centering
\caption{Comparison of out-of-sample prediction performance among five models}
\label{tab:Add-MINAR_out_of_sample_prediction}
\begin{tabular}{lccccc}
\toprule
Model & E1 & E2 & E3 & E4 & Number of Parameters \\
\midrule
MGINAR(1) & 1193.4011 & 23.3629 & 3.0965 & 0.8108 & 90 \\
iINAR(1) & 1078.6838 & 21.7809 & 2.2541 & 0.5324 & \textbf{18} \\
MINAR(1)\_ICLSE & 1074.5497 & 21.7362 & 2.4197 & 0.6049 & 27 \\
Add-MINAR(1)\_PROJ & 1073.6371 & 21.6774 & \textbf{1.6497} & \textbf{0.4742} & 36 \\
Add-MINAR(1)\_ICLSE & \textbf{1043.8693} & \textbf{21.2169} & 2.2032 & 0.5641 & 36 \\
\bottomrule
\end{tabular}
\end{table}

In out-of-sample prediction, the ranking of model performance undergoes notable changes, further highlighting critical differences in model generalization ability. As shown in Table \ref{tab:Add-MINAR_out_of_sample_prediction}, the Add-MINAR(1)\_ICLSE model demonstrates strong generalization performance, achieving the best results in both \( \mathrm{E}_1 \) (1043.8693) and \( \mathrm{E}_2 \) (21.2169). The Add-MINAR(1)\_PROJ model attains the optimal values in \( \mathrm{E}_3 \) (1.6497) and \( \mathrm{E}_4 \) (0.4742), although its \( \mathrm{E}_1 \) and \( \mathrm{E}_2 \) errors are slightly higher than those of the ICLSE counterpart. The structurally simple iINAR(1) model exhibits good generalization ability in terms of \( \mathrm{E}_3 \) and \( \mathrm{E}_4 \), even outperforming some more complex models. It is worth noting that the MINAR(1)\_ICLSE model, which performed excellently in-sample, shows a marked increase in certain metrics (e.g., \( \mathrm{E}_3 \) and \( \mathrm{E}_4 \)) in out-of-sample prediction, indicating diminished generalization capacity. The MGINAR(1) model, which has the largest number of parameters and the poorest in-sample performance, continues to exhibit the highest error metrics (e.g., \( \mathrm{E}_3 \) reaches 3.0965) in out-of-sample prediction, further confirming its overfitting problem. Synthesizing both in-sample and out-of-sample results, it can be seen that the Add-MINAR(1)-type models, particularly the one estimated via ICLSE, achieve a better overall balance between prediction accuracy and generalization robustness.

\begin{table}[H]
\centering
\caption{Coefficient matrix estimates of the Add-MINAR(1) model under two estimation methods}
\label{tab:Add-MINAR-para-transposed}
\begin{tabular}{lcc}
\toprule
Parameter Matrix & PROJ & ICLSE \\
\midrule
${\mathbf{A}}$ &
$\begin{array}{ccc}
0 & 0.0204 & 0.0267 \\
0.0538 & 0 & 0.1497 \\
0.2133 & 0.0649 & 0
\end{array}$ &
$\begin{array}{ccc}
0 & 0.0271 & 0.0334 \\
0.0512 & 0 & 0.1251 \\
0.2554 & 0.0837 & 0
\end{array}$ \\

\midrule

${\mathbf{B}}$ &
$\begin{array}{ccc}
0 & 0.1825 & 0.0738 \\
0.1400 & 0 & 0.1040 \\
0.1066 & 0.2081 & 0
\end{array}$ &
$\begin{array}{ccc}
0 & 0.2377 & 0.0785 \\
0.1385 & 0 & 0.0966 \\
0.0934 & 0.1925 & 0
\end{array}$ \\

\midrule

${\mathbf{C}}$ &
$\begin{array}{ccc}
0.0895 & 0.1418 & 0.1716 \\
0.5697 & 0.4829 & 0.4538 \\
0.5126 & 0.3322 & 0.3593\\
\end{array}$ &
$\begin{array}{ccc}
0.0947 & 0.2532 & 0.1609 \\
0.5989 & 0.5396 & 0.4550 \\
0.5013 & 0.4312 & 0.3774\\
\end{array}$ \\

\midrule

${\mathbf{D}}$ &
$\begin{array}{ccc}
3.3074 & 0.4192 & 3.5077 \\
2.5759 & 7.6996 & 15.3219 \\
6.2653 & 2.7847 & 26.3233\\
\end{array}$ &
$\begin{array}{ccc}
3.8677 & 1.6706 & 1.2176 \\
3.1119 & 1.4459 & 11.2514 \\
4.3406 & 3.2736 & 26.6276\\
\end{array}$ \\

\bottomrule
\end{tabular}
\end{table}

In the model specification, zero constraints are imposed on the diagonal elements of the coefficient matrices, i.e., \( \operatorname{diag}(\mathbf{A}) = \mathbf{0}_{m \times 1} \) and \( \operatorname{diag}(\mathbf{B}) = \mathbf{0}_{n \times 1} \), as model identification conditions. Table \ref{tab:Add-MINAR-para-transposed} presents the estimated parameter matrices of the Add-MINAR(1) model under the two estimation methods. It can be observed from the table that the row and column effect matrices \( \mathbf{A} \) and \( \mathbf{B} \) estimated by the ICLSE method for the Add-MINAR(1) model are generally higher than those estimated by PROJ, suggesting that the former tends to assign stronger explanatory weights to inter-variable relationships. This may reflect its stronger capacity to capture the persistence of system dynamics during the fitting process. Combining the in-sample and out-of-sample fitting results from Table \ref{tab:Add-MINAR_within_sample_fitting_new} and Table \ref{tab:Add-MINAR_out_of_sample_prediction}, it is also evident that the ICLSE method achieves superior in-sample fit.

\begin{figure}[H]
  \centering
  \includegraphics[width=1\textwidth]{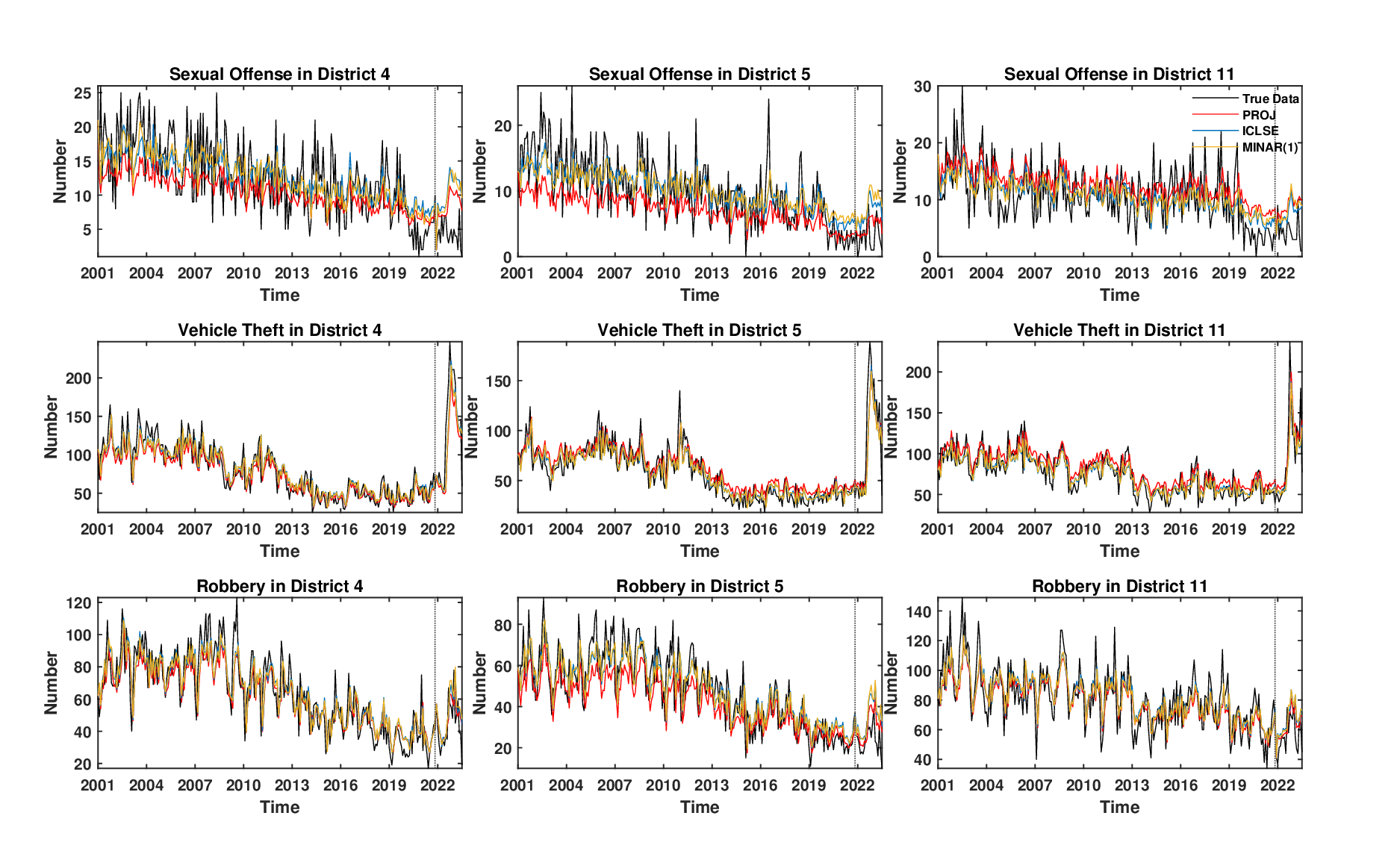}
  \caption{Comparison of observed data and predicted values from the Add-MINAR(1) model: The area to the right of the black dashed line indicates out-of-sample predictions.}
  \label{fig:Add-MINAR_prediction_plots}
\end{figure}

Figure \ref{fig:Add-MINAR_prediction_plots} presents the visualization of one-step-ahead predictions for the MGINAR(1), MINAR(1), and Add-MINAR(1) models, covering both in-sample and out-of-sample periods. The analysis reveals that, in predicting sexual offenses in Districts 4 and 5, the PROJ method for the Add-MINAR(1) model performs poorly, failing to effectively capture the data trends. The fitted values deviate substantially from the actual observations and are generally lower than the true values. Combined with the numerical results in Table \ref{tab:Add-MINAR_within_sample_fitting_new} and Table \ref{tab:Add-MINAR_out_of_sample_prediction}, it can be further confirmed that the predictive performance of Add-MINAR(1)\_PROJ is relatively low compared to the other models. In contrast, for vehicle theft and robbery in all three districts, all models exhibit reasonably good fitting ability and are able to capture the underlying trends of the data. Overall, the ICLSE version of the Add-MINAR(1) model demonstrates superior predictive accuracy for most crime types and districts, reflecting its stable advantage in modeling this type of crime data

To further intuitively illustrate the in-sample fitting superiority of the ICLSE method over the PROJ method, we evaluate the cumulative absolute error (CAE) over time.
The CAE is calculated by accumulating the absolute residuals at each time step; therefore, a lower curve indicates a smaller total fitting error over the sample period. 

 the blue curves (ICLSE) consistently remain below or perfectly overlap with the red curves (PROJ) across all panels. Notably, in the cases of Sexual Offenses across all districts, as well as Vehicle Theft and Robbery in District 5, the ICLSE method demonstrates a significantly slower error accumulation rate. This visual evidence perfectly corroborates the numerical results presented in Table \ref{tab:Add-MINAR_within_sample_fitting_new}, visually confirming that the ICLSE approach significantly reduces residual variance and yields more stable, accurate parameter estimates for the Add-MINAR(1) model compared to the PROJ method.


\section{Conclusion}
With the increasing complexity of data collection technologies and application scenarios, multi-dimensional integer time series frequently appear in numerous real-world applications such as public safety, transportation, and business analysis. Such data is generally presented in matrix form, for example, the sales panel data of "product types $\times$ sales outlets", which has a high degree of discrete data dimensions and contains a bidirectional dependency structure between rows and columns. Traditional vector modeling methods, when dealing with such data, not only destroy the original matrix structure, lose important row-column interaction information, but also face the "dimensionality disaster" where the number of parameters grows quadratically with the dimension, resulting in low parameter estimation efficiency, overfitting, and other problems. The MINAR model, although capable of directly handling matrix-form data and effectively depicting the relationship between rows and columns, has an undifferentiated bilinear product structure that cannot clearly separate the row and column effects in practice, and the model's interpretability is weak. To address these challenges, this paper expands the MINAR model based on the "additive structure" and provides a complete theoretical analysis and empirical research results. To address the problem of mixed row and column effects in the bilinear structure form and the insufficient model interpretability, this chapter proposes the Add-MINAR model. This model uses the additive structure to separate and model the row, column effects, and lag effects, and introduces a zero-diagonal constraint to ensure the unique identifiability of model parameters. This chapter estimates the parameters using the projection estimation method and the iterative conditional least squares method, and verifies the consistency and asymptotic normality of the two estimators. Simulation results show that the iterative conditional least squares estimation method has better estimation performance compared to the projection estimation method, and as the sample size increases, the estimation errors of the two methods tend to be consistent. In the Chicago real crime data, the Add-MINAR(1) model outperforms the iINAR(1), MGINAR(1), and MINAR(1) models in both in-sample fitting and out-of-sample prediction, especially in data with clear row-column structure, demonstrating strong explanatory ability and prediction stability.

\newpage
\bibliographystyle{elsarticle-num-names}
\bibliography{Add-MINAR_EDIT.bbl}
\newpage

\end{document}